# CURRENT STATUS DATA WITH COMPETING RISKS: CONSISTENCY AND RATES OF CONVERGENCE OF THE MLE

By Piet Groeneboom, Marloes H. Maathuis[1]
and Jon A. Wellner[2]

*Delft University of Technology and Vrije Universiteit Amsterdam, University of Washington and University of Washington*

We study nonparametric estimation of the sub-distribution functions for current status data with competing risks. Our main interest is in the nonparametric maximum likelihood estimator (MLE), and for comparison we also consider a simpler "naive estimator." Both types of estimators were studied by Jewell, van der Laan and Henneman [*Biometrika* (2003) **90** 183–197], but little was known about their large sample properties. We have started to fill this gap, by proving that the estimators are consistent and converge globally and locally at rate $n^{1/3}$. We also show that this local rate of convergence is optimal in a minimax sense. The proof of the local rate of convergence of the MLE uses new methods, and relies on a rate result for the sum of the MLEs of the sub-distribution functions which holds uniformly on a fixed neighborhood of a point. Our results are used in Groeneboom, Maathuis and Wellner [*Ann. Statist.* (2008) **36** 1064–1089] to obtain the local limiting distributions of the estimators.

**1. Introduction.** We study current status data with competing risks. Such data arise naturally in cross-sectional studies with several failure causes. Moreover, generalizations of these data arise in HIV vaccine trials (see [5]). The general framework is as follows. We analyze a system that can fail from $K$ competing risks, where $K \in \mathbb{N}$ is fixed. The random variables of interest are $(X, Y)$, where $X \in \mathbb{R}$ is the failure time of the system, and $Y \in \{1, \ldots, K\}$ is the corresponding failure cause. We cannot observe $(X, Y)$ directly. Rather, we observe the "current status" of the system at a single

Received September 2006; revised April 2007.
[1]Supported in part by NSF Grant DMS-02-03320.
[2]Supported in part by NSF Grants DMS-02-03320 and DMS-05-03822 and by NI-AID Grant 2R01 AI291968-04.

random time $T \in \mathbb{R}$, where $T$ is independent of $(X, Y)$. This means that at time $T$, we observe whether or not failure occurred, and if and only if failure occurred, we also observe the failure cause $Y$.

We want to estimate the bivariate distribution of $(X, Y)$. Since $Y \in \{1, \ldots, K\}$, this is equivalent to estimating the sub-distribution functions $F_{0k}(s) = P(X \leq s, Y = k)$, $k = 1, \ldots, K$. Note that the sum of the sub-distribution functions $\sum_{k=1}^{K} F_{0k}(s) = P(X \leq s)$ is the overall failure time distribution. This shows that the sub-distribution functions are related to each other and should be considered as a system.

We consider nonparametric estimation of the sub-distribution functions. This problem, or close variants thereof, has been studied by [5, 6, 7]. These papers introduced various nonparametric estimators, including the MLE (see [5, 7]) and a "naive estimator" (see [7]). They also provided algorithms to compute the estimators, and showed simulation studies that compared them. However, until now, little was known about the large sample properties of the estimators.

We have started to fill this gap by developing the local asymptotic theory for the MLE and the naive estimator. We study the MLE because it is a natural estimator that often exhibits good behavior. The simpler naive estimator was suggested to be asymptotically efficient for the estimation of smooth functionals [7], and we therefore consider it for comparison. In the present paper we prove consistency and rates of convergence. These results are used in [3] to obtain the local limiting distributions.

The outline of this paper is as follows. In Section 2 we introduce the estimators. We discuss their definitions, give existence and uniqueness results, and provide various characterizations in terms of necessary and sufficient conditions. Such characterizations are important since there is no closed form available for the MLE. In Section 3 we show that the estimators are globally and locally consistent. In Section 4 we prove that their global and local rates of convergence are $n^{1/3}$ (Theorems 4.1 and 4.17). We also prove that $n^{1/3}$ is an asymptotic local minimax lower bound for the rate of convergence (Proposition 4.4). Hence, the estimators converge locally at the optimal rate, in a minimax sense. The proof of the local rate of convergence of the MLE uses new methods. One of the main difficulties in this proof consists of handling the system of sub-distribution functions. We solve this problem by first deriving a rate result for the sum of the MLEs of the sub-distribution functions (Theorem 4.10). This rate result is stronger than usual, since it holds uniformly on a fixed neighborhood of a point, instead of on a shrinking neighborhood of order $n^{-1/3}$ (see Remark 4.11). Such a strong result is needed to handle potential sparsity of the jump points of the MLEs of the sub-distribution functions (see Remark 4.18). Technical proofs are collected in Section 5, and computational aspects of the estimators are discussed in the companion paper [3], Section 4.

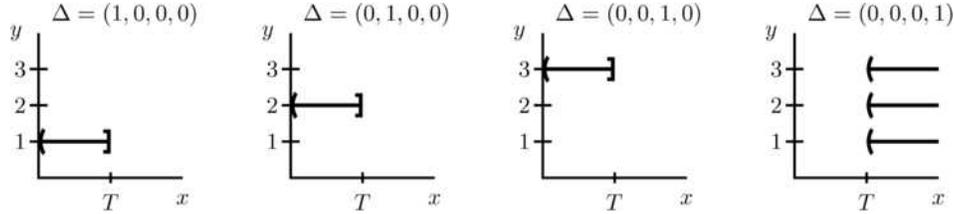

FIG. 1. *Graphical representation of the observed data $(T, \Delta)$ in an example with $K = 3$ competing risks. The black sets indicate the values of $(X, Y)$ that are consistent with $(T, \Delta)$, for each of the four possible values of $\Delta$.*

**2. The estimators.** We make the following assumptions: (a) the observation time $T$ is independent of the variables of interest $(X, Y)$, and (b) the system cannot fail from two or more causes at the same time. Assumption (a) is essential for the development of the theory. Assumption (b) ensures that the failure cause is well defined. This assumption is always satisfied by defining simultaneous failure from several causes as a new failure cause. We allow ties in the observation times.

We now introduce some notation. We denote the observed data by $(T, \Delta)$, where $T$ is the observation time and $\Delta = (\Delta_1, \ldots, \Delta_{K+1})$ is an indicator vector defined by $\Delta_k = 1\{X \leq T, Y = k\}$ for $k = 1, \ldots, K$, and $\Delta_{K+1} = 1\{X > T\}$. The observed data are illustrated in Figure 1. Let $(T_i, \Delta^i)$, $i = 1, \ldots, n$, be $n$ i.i.d. observations of $(T, \Delta)$, where $\Delta^i = (\Delta^i_1, \ldots, \Delta^i_{K+1})$. Note that we use the superscript $i$ as the index of an observation, and not as a power. The order statistics of $T_1, \ldots, T_n$ are denoted by $T_{(1)}, \ldots, T_{(n)}$. Furthermore, $G$ is the distribution of $T$, $G_n$ is the empirical distribution of $T_i$, $i, \ldots, n$, and $\mathbb{P}_n$ is the empirical distribution of $(T_i, \Delta^i)$, $i = 1, \ldots, n$. For any vector $(x_1, \ldots, x_K) \in \mathbb{R}^K$ we use the shorthand notation $x_+ = \sum_{k=1}^K x_k$, so that, for example, $\Delta_+ = \sum_{k=1}^K \Delta_k$ and $F_{0+}(s) = \sum_{k=1}^K F_{0k}(s)$. For any $K$-tuple $F = (F_1, \ldots, F_K)$ of sub-distribution functions, we define $F_{K+1}(s) = \int_{u>s} dF_+(u) = F_+(\infty) - F_+(s)$. Finally, we use the following conventions for indicator functions and integrals:

DEFINITION 2.1. Let $dA$ be a Lebesgue–Stieltjes measure. Then we define for $t < t_0$:

$$1_{[t_0, t)}(u) = -1_{[t, t_0)}(u) \quad \text{and} \quad \int_{[t_0, t)} f(u) \, dA(u) = -\int_{[t, t_0)} f(u) \, dA(u).$$

2.1. *Definitions of the estimators.* We first consider the MLE. To understand its form, let $F = (F_1, \ldots, F_K) \in \mathcal{F}_K$, where $\mathcal{F}_K$ is the collection of $K$-tuples $F = (F_1, \ldots, F_K)$ of sub-distribution functions on $\mathbb{R}$ with $F_+ \leq 1$. Under $F$ we have $\Delta | T \sim \text{Mult}_{K+1}(1, (F_1(T), \ldots, F_{K+1}(T)))$, so that the density of a single observation is given by

$$p_F(t, \delta) = \prod_{k=1}^{K+1} F_k(t)^{\delta_k} = \prod_{k=1}^{K} F_k(t)^{\delta_k} (1 - F_+(t))^{1-\delta_+}, \tag{1}$$

with respect to the dominating measure $\mu = G \times \#$, where $\#$ is the counting measure on $\{e_k : k = 1, \ldots, K+1\}$ and $e_k$ is the $k$th unit vector in $\mathbb{R}^{K+1}$. Hence, the log likelihood $l_n(F) = \int \log p_F(t, \delta) \, d\mathbb{P}_n(t, \delta)$ is given by

$$l_n(F) = \int \left\{ \sum_{k=1}^{K} \delta_k \log F_k(t) + (1 - \delta_+) \log(1 - F_+(t)) \right\} d\mathbb{P}_n(t, \delta). \tag{2}$$

It then follows that the MLE $\widehat{F}_n = (\widehat{F}_{n1}, \ldots, \widehat{F}_{nK})$ is defined by

$$l_n(\widehat{F}_n) = \max_{F \in \mathcal{F}_K} l_n(F). \tag{3}$$

The naive estimator $\widetilde{F}_n = (\widetilde{F}_{n1}, \ldots, \widetilde{F}_{nK})$ is defined by

$$l_{nk}(\widetilde{F}_{nk}) = \max_{F_k \in \mathcal{F}} l_{nk}(F_k), \qquad k = 1, \ldots, K, \tag{4}$$

where $\mathcal{F}$ is the collection of all distribution functions on $\mathbb{R}$ and $l_{nk}(\cdot)$ is the marginal log likelihood for the reduced current status data $(T_i, \Delta_k^i)$, $i = 1, \ldots, n$:

$$l_{nk}(F_k) = \int \{\delta_k \log F_k(t) + (1 - \delta_k) \log(1 - F_k(t))\} \, d\mathbb{P}_n(t, \delta),$$

$$k = 1, \ldots, K.$$

Thus, $\widetilde{F}_{nk}$ uses only the $k$th entry of the $\Delta$-vector. We see that the naive estimator splits the estimation problem into $K$ well-known univariate current status problems. Therefore, its computation and asymptotic theory follow straightforwardly from known results on current status data. But this simplification comes at a cost. For example, it follows immediately that the constraint $\widetilde{F}_{n+} \leq 1$ may be violated (see [7]).

We note that both $\widehat{F}_{n+}$ and $\widetilde{F}_{n+}$ provide estimators for the overall failure time distribution $F_{0+}$. A third estimator for this distribution is given by the MLE for the reduced current status data $(T, \Delta_+)$, ignoring information on the failure causes. These three estimators are typically not the same (see [5]).

To compare the MLE and the naive estimator, we now define the naive estimator by a single optimization problem:

$$\widetilde{l}_n(\widetilde{F}_n) = \max_{F \in \mathcal{F}^K} \widetilde{l}_n(F) \qquad \text{where } \widetilde{l}_n(F) = \sum_{k=1}^{K} l_{nk}(F_k),$$

and $\mathcal{F}^K$ is the $K$-fold product of $\mathcal{F}$. By comparing this to the optimization problem for the MLE, we note the following differences:

(a) The object function $l_n(F)$ for the MLE contains the term $1 - F_+$, involving the sum of the sub-distribution functions, while the object function $\widetilde{l}_n(F)$ for the naive estimator only contains the individual components.

(b) The space $\mathcal{F}_K$ for the MLE contains the constraint $F_+ \leq 1$, while the space $\mathcal{F}^K$ for the naive estimator only involves the individual components.

The more complicated object function for the MLE forces us to work with the system of sub-distribution functions, and poses new challenges in the derivation of the local rate of convergence of the MLE. Moreover, it gives rise to a new self-induced limiting process for the local limiting distribution of the MLE (see [3]). The constraint $F_+ \leq 1$ on the space over which we maximize is important for small sample sizes, but its effect vanishes asymptotically. These observations are supported by simulations in [3], Section 4.

2.2. *Existence and uniqueness.* Since only values of the sub-distribution functions at the observation times appear in the log likelihoods $l_{nk}(F_k)$ and $l_n(F)$, we limit ourselves to estimating these values. This means that the optimization problems (3) and (4) reduce to finite-dimensional optimization problems. Hence, their solutions exist by [19], Corollary 38.10.

For the naive estimator, the values of the sub-distribution functions at all observation times enter in the log likelihood $l_{nk}(F_k)$. Together with strict concavity of $l_{nk}(F_k)$, this implies that $\widetilde{F}_{nk}$ is unique at all observation times, for $k = 1, \ldots, K$. For the MLE, $F_k(T_i)$ appears in the log likelihood $l_n(F)$ if and only if $\Delta_k^i + \Delta_{K+1}^i > 0$. This motivates the following definition and result:

DEFINITION 2.2. For each $k = 1, \ldots, K+1$, we define the set $\mathcal{T}_k$ by

(5) $$\mathcal{T}_k = \{T_i, i = 1, \ldots, n : \Delta_k^i + \Delta_{K+1}^i > 0\} \cup \{T_{(n)}\}.$$

PROPOSITION 2.3. *For each $k = 1, \ldots, K+1$, $\widehat{F}_{nk}(t)$ is unique at $t \in \mathcal{T}_k$. Moreover, $\widehat{F}_{nk}(\infty)$ is unique if and only if $\Delta_{K+1}^i = 0$ for all observations with $T_i = T_{(n)}$.*

PROOF. We first prove uniqueness of $\widehat{F}_{nk}(t)$ at $t \in \mathcal{T}_k$, for $k = 1, \ldots, K$. Let $k \in \{1, \ldots, K\}$. Strict concavity of the log likelihood immediately gives uniqueness of $\widehat{F}_{nk}$ at points $T_i$ with $\Delta_k^i = 1$. Note that the log likelihood is not *strictly* concave in $\widehat{F}_{nk}(T_i)$ if $\Delta_{K+1}^i = 1$, so that we need to do more work to prove uniqueness at these points. First, one can show that $\widehat{F}_{nk}$ can only assign mass to intervals of the following form:

(i) $(T_i, T_j]$ where $\Delta_{K+1}^i = 1$, $\Delta_k^j = 1$ and $\Delta_k^\ell = \Delta_{K+1}^\ell = 0$ for all $\ell$ such that $T_i < T_\ell < T_j$,

(ii) $(T_i, \infty)$ where $T_i = T_{(n)}$ and $\Delta^i_{K+1} = 1$

(see [5], Lemma 1, or use the concept of the height map of [9]). Note that $\widehat{F}_{nk}$ is unique at the right endpoints of the intervals given in (i), since $\widehat{F}_{nk}$ is unique at points $T_i$ with $\Delta^i_k = 1$. This implies that the probability mass in each interval given in (i) is unique. In turn, this implies that $\widehat{F}_{nk}$ is unique at all points that are not in the interior of these intervals. In particular, this gives uniqueness of $\widehat{F}_{nk}(t)$ at $t \in \mathcal{T}_k$. The uniqueness statement about $\widehat{F}_{n,K+1}$ follows from the uniqueness of $\widehat{F}_{n1}, \ldots, \widehat{F}_{nK}$.

We now prove the statement about $\widehat{F}_{nk}(\infty)$. First, if $\Delta^i_{K+1} = 0$ for all observations with $T_i = T_{(n)}$, then $\widehat{F}_{nk}$ can only assign mass to the intervals given in (i). Hence, $\widehat{F}_{nk}(\infty) = \widehat{F}_{nk}(T_{(n)})$, and since $\widehat{F}_{nk}(T_{(n)})$ was already proved to be unique, it follows that $\widehat{F}_{nk}(\infty)$ is unique. Conversely, if there is a $T_i = T_{(n)}$ with $\Delta^i_{K+1} = 1$, then the log likelihood contains the term $\log(1 - F_+(T_{(n)}))$. Hence, $\widehat{F}_{n+}$ must assign mass to the right of $T_{(n)}$ in order to get $l_n(\widehat{F}_n) > -\infty$. The MLE is indifferent to the distribution of this mass over $\widehat{F}_{n1}, \ldots, \widehat{F}_{nK}$, since their separate contributions do not appear in the log likelihood. Hence, $\widehat{F}_{nk}(\infty)$ is nonunique in this case. □

2.3. *Characterizations.* Characterizations of the naive estimators $\widetilde{F}_{n1}, \ldots, \widetilde{F}_{nK}$ follow from [4], Propositions 1.1 and 1.2, pages 39–41. Characterizations of the MLE can be derived from Karush–Kuhn–Tucker conditions, since the optimization problem can be reduced to a finite-dimensional optimization problem (see the first paragraph of Section 2.2). However, we give characterizations with direct proofs. These methods do not use the discrete nature of the problem, so that they can also be used for truly infinite-dimensional optimization problems.

DEFINITION 2.4. We define the processes $V_{nk}$ by

$$(6) \qquad V_{nk}(t) = \int_{u \leq t} \delta_k \, d\mathbb{P}_n(u, \delta), \qquad t \in \mathbb{R}, \ k = 1, \ldots, K+1.$$

Moreover, let $\bar{\mathcal{F}}_K$ be the collection of $K$-tuples of bounded nonnegative nondecreasing right-continuous functions.

Using this notation, we can write $l_n(F) = \sum_{k=1}^{K+1} \int \log F_k(u) \, dV_{nk}(u)$. In Lemma 2.5 we translate the optimization problem (3) into an optimization problem over a cone, by removing the constraint $F_+ \leq 1$. Subsequently, we give a basic characterization in Proposition 2.6. This characterization leads to various corollaries, of which Corollary 2.10 is most important for the sequel.

LEMMA 2.5. *$\widehat{F}_n$ maximizes $l_n(F)$ over $\mathcal{F}_K$ if and only if $\widehat{F}_n$ maximizes $\bar{l}_n(F)$ over $\bar{\mathcal{F}}_K$, where*

$$\bar{l}_n(F) = \sum_{k=1}^{K+1} \int \log F_k(u) \, dV_{nk}(u) - F_+(\infty).$$

PROOF. (Necessity.) Let $\widehat{F}_n$ maximize $l_n(F)$ over $\mathcal{F}_K$, and let $F \in \bar{\mathcal{F}}_K$. We want to show that $\bar{l}_n(\widehat{F}_n) \geq \bar{l}_n(F)$. Note that this inequality holds trivially if $F_+(\infty) = 0$. Hence, we assume $F_+(\infty) = c > 0$. Then $F/c \in \mathcal{F}_K$, and $l_n(\widehat{F}_n) \geq l_n(F/c)$, by the assumption that $\widehat{F}_n$ maximizes $l_n(F)$ over $\mathcal{F}_K$. Together with $\widehat{F}_{n+}(\infty) = 1$ this yields

$$\bar{l}_n(\widehat{F}_n) = l_n(\widehat{F}_n) - 1 \geq l_n(F/c) - 1$$
$$= \sum_{k=1}^{K+1} \int \log F_k(u) \, dV_{nk}(u) - \log c - 1$$
$$= \bar{l}_n(F) + c - \log c - 1 \geq \bar{l}_n(F).$$

The last inequality follows since $x - \log x - 1 \geq 0$ for $x > 0$.

(Sufficiency.) Let $\widehat{F}_n$ maximize $\bar{l}_n(F)$ over $\bar{\mathcal{F}}_K$, and let $\widehat{F}_{n+}(\infty) = c$. As before, we may assume $c > 0$. Then $\bar{l}_n(\widehat{F}_n) \geq \bar{l}_n(\widehat{F}_n/c)$, and by the same reasoning as above this gives $\bar{l}_n(\widehat{F}_n) \geq \bar{l}_n(\widehat{F}_n/c) = l_n(\widehat{F}_n/c) - 1 = \bar{l}_n(\widehat{F}_n) + c - \log c - 1$. Since $x - \log x - 1 \leq 0$ if and only if $x = 1$, this yields $c = 1$. Hence, $\widehat{F}_n \in \mathcal{F}_K$, and $\widehat{F}_n$ maximizes $l_n(F)$ over $\mathcal{F}_K \subset \bar{\mathcal{F}}_K$. □

We now obtain the following basic characterization of the MLE.

PROPOSITION 2.6. *$\widehat{F}_n$ maximizes $l_n(F)$ over $\mathcal{F}_K$ if and only if $\widehat{F}_n \in \bar{\mathcal{F}}_K$ and the following two conditions hold for all $k = 1, \ldots, K$:*

$$(7) \qquad \int_{u \geq t} \frac{dV_{nk}(u)}{\widehat{F}_{nk}(u)} + \int_{u < t} \frac{dV_{n,K+1}(u)}{\widehat{F}_{n,K+1}(u)} \leq 1, \qquad t \in \mathbb{R},$$

$$(8) \qquad \int \left\{ \int_{u \geq t} \frac{dV_{nk}(u)}{\widehat{F}_{nk}(u)} + \int_{u < t} \frac{dV_{n,K+1}(u)}{\widehat{F}_{n,K+1}(u)} - 1 \right\} d\widehat{F}_{nk}(t) = 0.$$

PROOF. (Necessity.) Let $\widehat{F}_n$ maximize $l_n(F)$ over $\mathcal{F}_K$. Then $\widehat{F}_n$ also maximizes $\bar{l}_n(F)$ over $\bar{\mathcal{F}}_K$, by Lemma 2.5. Fix $k \in \{1, \ldots, K\}$, and define the perturbation $\widehat{F}_n^{(h)} = (\widehat{F}_{n1}^{(h)}, \ldots, \widehat{F}_{nK}^{(h)})$ by $\widehat{F}_{nk}^{(h)} = (1+h)\widehat{F}_{nk}$ and $\widehat{F}_{nj}^{(h)} = \widehat{F}_{nj}$ for $j \neq k$. Since $\widehat{F}_n^{(h)} \in \bar{\mathcal{F}}_K$ for $|h| < 1$, we get

$$0 = \lim_{h \to 0} h^{-1}\{\bar{l}_n(\widehat{F}_n^{(h)}) - \bar{l}_n(\widehat{F}_n)\}$$

$$= \int dV_{nk}(u) + \int \frac{\widehat{F}_{nk}(\infty) - \widehat{F}_{nk}(u)}{\widehat{F}_{n,K+1}(u)} dV_{n,K+1}(u) - \widehat{F}_{nk}(\infty)$$

$$= \int \left\{ \int_{u \geq t} \frac{dV_{nk}(u)}{\widehat{F}_{nk}(u)} + \int_{u < t} \frac{dV_{n,K+1}(u)}{\widehat{F}_{n,K+1}(u)} - 1 \right\} d\widehat{F}_{nk}(t),$$

using Fubini's theorem to obtain the last line. This gives condition (8). Next, let $t \in \mathbb{R}$, and define the perturbation $\widehat{F}_n^{(h,t)} = (\widehat{F}_{n1}^{(h,t)}, \ldots, \widehat{F}_{nK}^{(h,t)})$ by $\widehat{F}_{nk}^{(h,t)}(u) = \widehat{F}_{nk}(u) + h1_{[t,\infty)}(u)$ and $\widehat{F}_{nj}^{(h,t)} = \widehat{F}_{nj}$ for $j \neq k$. Since $\widehat{F}_n^{(h,t)} \in \bar{\mathcal{F}}_K$ for $h \geq 0$, we get

$$0 \geq \lim_{h \downarrow 0} h^{-1} \{ \bar{l}_n(\widehat{F}_n^{(h,t)}) - \bar{l}_n(\widehat{F}_n) \}$$

$$= \int_{u \geq t} \frac{dV_{nk}(u)}{\widehat{F}_{nk}(u)} + \int_{u < t} \frac{dV_{n,K+1}(u)}{\widehat{F}_{n,K+1}(u)} - 1,$$

which is condition (7).

(Sufficiency.) Let $\widehat{F}_n \in \bar{\mathcal{F}}_K$ satisfy conditions (7) and (8), and let $F \in \bar{\mathcal{F}}_K$. We want to show that $\bar{l}_n(\widehat{F}_n) \geq \bar{l}_n(F)$. Concavity of the logarithm yields

$$\bar{l}_n(F) - \bar{l}_n(\widehat{F}_n) \leq \sum_{k=1}^{K+1} \int \frac{F_k(u) - \widehat{F}_{nk}(u)}{\widehat{F}_{nk}(u)} dV_{nk}(u) - F_+(\infty) + \widehat{F}_{n+}(\infty).$$

We now show that the right-hand side of this display is nonpositive. By Fubini, we have

$$\sum_{k=1}^{K} \int \frac{F_k(u) - \widehat{F}_{nk}(u)}{\widehat{F}_{nk}(u)} dV_{nk}(u) = \sum_{k=1}^{K} \int \int_{t \leq u} d(F_k - \widehat{F}_{nk})(t) \frac{dV_{nk}(u)}{\widehat{F}_{nk}(u)}$$

$$= \sum_{k=1}^{K} \int \int_{u \geq t} \frac{dV_{nk}(u)}{\widehat{F}_{nk}(u)} d(F_k - \widehat{F}_{nk})(t)$$

and

$$\int \frac{F_{K+1}(u) - \widehat{F}_{n,K+1}(u)}{\widehat{F}_{n,K+1}(u)} dV_{n,K+1}(u) = \int \int_{t > u} d(F_+ - \widehat{F}_{n+})(t) \frac{dV_{n,K+1}(u)}{\widehat{F}_{n,K+1}(u)}$$

$$= \sum_{k=1}^{K} \int \int_{u < t} \frac{dV_{n,K+1}(u)}{\widehat{F}_{n,K+1}(u)} d(F_k - \widehat{F}_{nk})(t).$$

Combining the last three displays gives

$$\bar{l}_n(F) - \bar{l}_n(\widehat{F}_n) \leq \sum_{k=1}^{K} \int \left\{ \int_{u \geq t} \frac{dV_{nk}(u)}{\widehat{F}_{nk}(u)} + \int_{u < t} \frac{dV_{n,K+1}(u)}{\widehat{F}_{n,K+1}(u)} - 1 \right\} d(F_k - \widehat{F}_{nk})(t)$$

$$= \sum_{k=1}^{K} \int \left\{ \int_{u \geq t} \frac{dV_{nk}(u)}{\widehat{F}_{nk}(u)} + \int_{u < t} \frac{dV_{n,K+1}(u)}{\widehat{F}_{n,K+1}(u)} - 1 \right\} dF_k(t) \leq 0,$$

where the equality follows from (8), and the final inequality follows from (7). Hence $\widehat{F}_n$ maximizes $\bar{l}_n(F)$ over $\bar{\mathcal{F}}_K$, and by Lemma 2.5 this implies that $\widehat{F}_n$ maximizes $l_n(F)$ over $\mathcal{F}_K$. □

DEFINITION 2.7. We say that $t$ is a *point of increase* of a right-continuous function $F$ if $F(t) > F(t-\varepsilon)$ for every $\varepsilon > 0$ (note that this definition is slightly different from the usual definition). Moreover, for $F \in \bar{\mathcal{F}}_K$, we define

$$(9) \qquad \beta_{nF} = 1 - \int \frac{dV_{n,K+1}(u)}{F_{K+1}(u)}.$$

Note that $\beta_{n\widehat{F}_n}$ is uniquely defined, since $\widehat{F}_{n,K+1}(t)$ is unique at points $t$ where $dV_{n,K+1}$ has mass (Proposition 2.3). We now rewrite the characterization in Proposition 2.6 in terms of $\beta_{n\widehat{F}_n}$:

COROLLARY 2.8. $\widehat{F}_n$ maximizes $l_n(F)$ over $\mathcal{F}_K$ if and only if $\widehat{F}_n \in \bar{\mathcal{F}}_K$ and the following holds for all $k = 1, \ldots, K$:

$$(10) \qquad \int_{u \geq t} \left\{ \frac{dV_{nk}(u)}{\widehat{F}_{nk}(u)} - \frac{dV_{n,K+1}(u)}{\widehat{F}_{n,K+1}(u)} \right\} \leq \beta_{n\widehat{F}_n}, \qquad t \in \mathbb{R},$$

where equality holds if $t$ is a point of increase of $\widehat{F}_{nk}$.

PROOF. Since the integrand of (8) is a left-continuous function of $t$, conditions (7) and (8) of Proposition 2.6 are equivalent to the condition that for all $k = 1, \ldots, K$,

$$\int_{u \geq t} \frac{dV_{nk}(u)}{\widehat{F}_{nk}(u)} + \int_{u < t} \frac{dV_{n,K+1}(u)}{\widehat{F}_{n,K+1}(u)} \leq 1, \qquad t \in \mathbb{R},$$

where equality must hold if $t$ is a point of increase of $\widehat{F}_{nk}$. Combining this with

$$\int_{u < t} \frac{dV_{n,K+1}(u)}{\widehat{F}_{n,K+1}(u)} = 1 - \beta_{n\widehat{F}_n} - \int_{u \geq t} \frac{dV_{n,K+1}(u)}{\widehat{F}_{n,K+1}(u)}, \qquad t \in \mathbb{R}$$

completes the proof. □

We determine the sign of $\beta_{n\widehat{F}_n}$ in Corollary 2.9:

COROLLARY 2.9. Let $\widehat{F}_n$ maximize $l_n(F)$ over $\mathcal{F}_K$. Then $\beta_{n\widehat{F}_n} \geq 0$, and $\beta_{n\widehat{F}_n} = 0$ if and only if there is an observation with $T_i = T_{(n)}$ and $\Delta^i_{K+1} = 1$.

PROOF. Taking $t > T_{(n)}$ in Corollary 2.8 implies that $\beta_{n\widehat{F}_n} \geq 0$. Now suppose that there is a $T_i = T_{(n)}$ with $\Delta^i_{K+1} = 1$. Then we must have $\widehat{F}_{n+}(T_{(n)}) < 1$ to obtain $l_n(\widehat{F}_n) > -\infty$. Hence, there must be a $k \in \{1, \ldots, K\}$ such that $\widehat{F}_{nk}$ has points of increase $t > T_{(n)}$. Corollary 2.8 then implies that $\beta_{n\widehat{F}_n} = 0$. Next, suppose that there does not exist a $T_i = T_{(n)}$ with $\Delta^i_{K+1} = 1$. Then

$$\int_{u \geq T_{(n)}} \left\{ \frac{dV_{nk}(u)}{\widehat{F}_{nk}(u)} - \frac{dV_{n,K+1}(u)}{\widehat{F}_{n,K+1}(u)} \right\} = \int_{u \geq T_{(n)}} \frac{dV_{nk}(u)}{\widehat{F}_{nk}(u)} > 0,$$

and by Corollary 2.8 this implies $\beta_{n\widehat{F}_n} > 0$. □

We now make a first step toward localizing the characterization, in Corollary 2.10. This corollary forms the basis of Proposition 4.8, which is used in the proofs of the local rate of convergence and the limiting distribution of the MLE.

COROLLARY 2.10. $\widehat{F}_n$ maximizes $l_n(F)$ over $\mathcal{F}_K$ if and only if $\widehat{F}_n \in \bar{\mathcal{F}}_K$ and the following holds for all $k = 1, \ldots, K$ and each point of increase $\tau_{nk}$ of $\widehat{F}_{nk}$:

$$(11) \quad \int_{[\tau_{nk}, s)} \left\{ \frac{dV_{nk}(u)}{\widehat{F}_{nk}(u)} - \frac{dV_{n,K+1}(u)}{\widehat{F}_{n,K+1}(u)} \right\} \geq \beta_{n\widehat{F}_n} 1_{[\tau_{nk}, s)}(T_{(n)}), \qquad s \in \mathbb{R},$$

where equality holds if $s$ is a point of increase of $\widehat{F}_{nk}$, and if $s > T_{(n)}$.

PROOF. Let $\widehat{F}_n$ maximize $l_n(\cdot)$ over $\mathcal{F}_K$. Let $s > \tau_{nk}$. If $\tau_{nk} < s \leq T_{(n)}$, then (11) follows by applying (10) to $t = \tau_{nk}$ and $t = s$, and subtracting the resulting equations. If $\tau_{nk} \leq T_{(n)} < s$, then

$$\int_{[\tau_{nk}, s)} \left\{ \frac{dV_{nk}(u)}{\widehat{F}_{nk}(u)} - \frac{dV_{n,K+1}(u)}{\widehat{F}_{n,K+1}(u)} \right\} = \int_{u \geq \tau_{nk}} \left\{ \frac{dV_{nk}(u)}{\widehat{F}_{nk}(u)} - \frac{dV_{n,K+1}(u)}{\widehat{F}_{n,K+1}(u)} \right\},$$

so that the statement follows by applying (10) to $t = \tau_{nk}$. If $T_{(n)} < \tau_{nk} < s$, then the left-hand side of (10) equals zero for $t = \tau_{nk}$ and $t = s$. The inequalities for $s < \tau_{nk}$ can be derived analogously. Finally, the inequality (11) and the corresponding equality condition imply (10). □

**3. Consistency.** Hellinger and $L_r(G)$ $(r \geq 1)$ consistency of the naive estimator follow from [13, 18]. Local consistency of the naive estimator follows from [4, 13]. In this section we prove similar results for the MLE. First, note that for two vectors of functions $F = (F_1, \ldots, F_K)$ and $F_0 = (F_{01}, \ldots, F_{0K})$ in $\mathcal{F}_K$, the Hellinger distance $h(p_F, p_{F_0})$ and the total variation distance $d_{\mathrm{TV}}(p_F, p_{F_0})$ in our model are given by

$$(12) \quad h^2(p_F, p_{F_0}) = \tfrac{1}{2} \int (\sqrt{p_F} - \sqrt{p_{F_0}})^2 \, d\mu = \tfrac{1}{2} \sum_{k=1}^{K+1} \int (\sqrt{F_k} - \sqrt{F_{0k}})^2 \, dG,$$

$$(13) \quad d_{\mathrm{TV}}(p_F, p_{F_0}) = \tfrac{1}{2} \sum_{k=1}^{K+1} \int |F_k - F_{0k}| \, dG,$$

where $\mu = G \times \#$, and $p_F$ and $\#$ are defined in (1). The MLE is Hellinger consistent:

THEOREM 3.1. $h(p_{\widehat{F}_n}, p_{F_0}) \to_{a.s.} 0$.

PROOF. Since $\mathcal{P} = \{p_F : F \in \mathcal{F}_K\}$ is convex, we can use the following inequality:

$$h^2(p_{\widehat{F}_n}, p_{F_0}) \le (\mathbb{P}_n - P) \phi(p_{\widehat{F}_n}/p_{F_0}),$$

where $\phi(t) = (t-1)/(t+1)$ ([18], Proposition 3; see also [11] and [14, 15]). Hence, it is sufficient to prove that $\{\phi(p_F/p_{F_0}) : F \in \mathcal{F}_K\}$ is a $P$-Glivenko–Cantelli class. This can be shown by Glivenko–Cantelli preservation theorems of [18], using indicators of $VC$-classes of sets and monotone functions as building blocks. Alternatively, the result follows directly from [18], Theorem 9 by viewing the problem as a bivariate censored data problem for $(X, Y)$. □

$L_r(G)$ consistency is given in Corollary 3.2, where the $L_r(G)$ distance is defined by

$$(14) \quad \|F - F_0\|_{G,r}^r = \sum_{k=1}^{K+1} \int |F_k(t) - F_{0k}(t)|^r \, dG(t), \qquad r \ge 1.$$

COROLLARY 3.2. $\|\widehat{F}_n - F_0\|_{G,r} \to_{a.s.} 0$ for $r \ge 1$.

PROOF. Note that $\|F - F_0\|_{G,1} = 2 d_{\mathrm{TV}}(p_F, p_{F_0})$. Hence, the statement for $r = 1$ follows from the well-known inequality $d_{\mathrm{TV}}(p_{F_1}, p_{F_2}) \le \sqrt{2} h(p_{F_1}, p_{F_2})$. The result for $r > 1$ follows from $|a - b|^r \le |a - b|$ for $a, b \in [0, 1]$ and $r > 1$. □

Note that Theorem 3.1 and Corollary 3.2 hold without any additional assumptions. The quantities in these statements are integrated with respect to $G$, showing the importance of the observation time distribution. For example, the results do not imply consistency at intervals where $G$ has zero mass. Such issues should be taken into account if $G$ can be chosen by design.

Under some additional assumptions, Maathuis ([10], Section 4.2) proved several forms of local and uniform consistency using methods from [13], Section 3. One such result is needed in the proof of the local rate of convergence of the MLE, and is given below:

PROPOSITION 3.3. *Let $F_{01}, \ldots, F_{0K}$ be continuous at $t_0$, and let $G$ be continuously differentiable at $t_0$ with strictly positive derivative $g(t_0)$. Then there exists an $r > 0$ such that*

$$\sup_{t \in [t_0-r, t_0+r]} |\widehat{F}_{nk}(t) - F_{0k}(t)| \to_{a.s.} 0, \qquad k = 1, \ldots, K.$$

PROOF. Let $k \in \{1, \ldots, K\}$ and choose the constant $r > 0$ such that $F_{0k}$ is continuous on $[t_0 - 2r, t_0 + 2r]$ and $g(t) > g(t_0)/2$ for $t \in [t_0 - 2r, t_0 + 2r]$. Fix an $\omega$ for which the $L_1(G)$ consistency holds, and suppose there is an $x_0 \in [t_0 - r, t_0 + r]$ for which $\widehat{F}_{nk}(x_0, \omega)$ does not converge to $F_{0k}(x_0)$. Then there is an $\varepsilon > 0$ such that for all $n_1 > 0$ there is an $n > n_1$ such that $|\widehat{F}_{nk}(x_0, \omega) - F_{0k}(x_0)| > \varepsilon$. Using the monotonicity of $\widehat{F}_{nk}$ and the continuity of $F_{0k}$, this implies there is a $\gamma > 0$ such that $|\widehat{F}_{nk}(t, \omega) - F_{0k}(t)| > \varepsilon/2$ for all $t \in (x_0 - \gamma, x_0]$ or $[x_0, x_0 + \gamma)$ and $[x_0 - \gamma, x_0 + \gamma] \subset [t_0 - 2r, t_0 + 2r]$. This yields that $\int |\widehat{F}_{nk}(t, \omega) - F_{0k}(t)| \, dG(t) > \gamma \varepsilon g(t_0)/4$, which contradicts $L_1(G)$ consistency. Uniform consistency follows since $F_{0k}$ is continuous. □

**4. Rate of convergence.** The Hellinger rate of convergence of the naive estimator is $n^{1/3}$. This follows from [15] or [17], Theorem 3.4.4, page 327. Under certain regularity conditions, the local rate of convergence of the naive estimator is also $n^{1/3}$; see [4], Lemma 5.4, page 95. This local rate result implies that the distance between two successive jump points of $\widetilde{F}_{nk}$ around a point $t_0$ is of order $O_p(n^{-1/3})$.

In this section we discuss similar results for the MLE. In Section 4.1 we show that the global rate of convergence is $n^{1/3}$. In Section 4.2 we prove that $n^{1/3}$ is an asymptotic local minimax lower bound for the rate of convergence, meaning that no estimator can converge locally at a rate faster than $n^{1/3}$, in a minimax sense. Hence, the naive estimator converges locally at the optimal rate. Since the MLE is expected to be at least as good as the naive estimator, one may expect that the MLE also converges locally at the optimal rate of $n^{1/3}$. This is indeed the case, and this is proved in Section 4.3 (Theorem 4.17). Our main tool for proving this result is Theorem 4.10, which gives a uniform rate of convergence of $\widehat{F}_{n+}$ on a fixed neighborhood of a point, rather than on the usual shrinking neighborhood of order $n^{-1/3}$. Such a strong rate result is needed to handle potential sparsity of the jump points of the MLEs of the sub-distribution functions (see Remark 4.18). Some technical proofs are deferred to Section 5.

4.1. *Global rate of convergence.*

THEOREM 4.1. $n^{1/3} h(p_{\widehat{F}_n}, p_{F_0}) = O_p(1)$.

PROOF. We use the rate theorem of Van der Vaart and Wellner ([17], Theorem 3.4.1, page 322) with

$$m_{p_F}(t, \delta) = \log\left(\frac{p_F(t, \delta) + p_{F_0}(t, \delta)}{2 p_{F_0}(t, \delta)}\right),$$

$\mathbb{M}_n(F) = \mathbb{P}_n m_{p_F}$, $M(F) = P m_{p_F}$ and $\mathbb{G}_n m_{p_F} = \sqrt{n}(\mathbb{M}_n - M)(F)$. The key condition to verify is $E\|\mathbb{G}_n\|_{\mathcal{M}_\gamma} \lesssim \phi_n(\gamma)$, where $\mathcal{M}_\gamma = \{m_{p_F} - m_{p_{F_0}} : h(p_F, p_{F_0}) < \gamma\}$ and $\phi_n(\gamma)/\gamma^\alpha$ is a decreasing function in $\gamma$ for some $\alpha < 2$. For this purpose we use Theorem 3.4.4 of [17], which states that the functions $m_{p_F}$ fit the setup of Theorem 3.4.1 of [17], and that

(15) $\qquad E\|\mathbb{G}_n\|_{\mathcal{M}_\gamma} \leq \widetilde{J}_{[\,]}(\gamma, \mathcal{P}, h)\{1 + \widetilde{J}_{[\,]}(\gamma, \mathcal{P}, h)\gamma^{-2} n^{-1/2}\},$

where $\widetilde{J}_{[\,]}(\gamma, \mathcal{P}, h) = \int_0^\gamma \sqrt{1 + \log N_{[\,]}(\varepsilon, \mathcal{P}, h)} \, d\varepsilon$ and $\log N_{[\,]}(\varepsilon, \mathcal{P}, h)$ is the $\varepsilon$-entropy with bracketing for $\mathcal{P} = \{p_F : F \in \mathcal{F}_K\}$ with respect to Hellinger distance $h$. We first bound the bracketing number $N_{[\,]}(\varepsilon, \mathcal{P}, h)$. Let $F = (F_1, \ldots, F_K) \in \mathcal{F}_K$. For each $k = 1, \ldots, K+1$, let $[l_k, u_k]$ be a bracket containing $F_k$, with size $\int(\sqrt{u_k} - \sqrt{l_k})^2 \, dG \leq \varepsilon^2/(K+1)$. Then

$$[p_l(t, \delta), p_u(t, \delta)] = \left[\prod_{k=1}^{K+1} l_k(t)^{\delta_k}, \prod_{k=1}^{K+1} u_k(t)^{\delta_k}\right]$$

is a bracket containing $p_F$, and its Hellinger size is bounded by $\varepsilon$.

Note that all $F_k$, $k = 1, \ldots, K+1$, are contained in the class $\mathcal{F} = \{F : \mathbb{R} \mapsto [0, 1] \text{ is monotone}\}$, and it is well known that $\log N_{[\,]}(\delta, \mathcal{F}, L_2(Q)) \lesssim 1/\delta$, uniformly in $Q$. Hence, considering all possible combinations of $(K+1)$-tuples of the brackets $[l_k, u_k]$, it follows that

$$\log N_{[\,]}(\varepsilon, \mathcal{P}, h) \leq \log(\{N_{[\,]}(\varepsilon/\sqrt{K+1}, \mathcal{F}, L_2(G))\}^{K+1})$$
$$= (K+1) \log N_{[\,]}(\varepsilon/\sqrt{K+1}, \mathcal{F}, L_2(G)) \lesssim (K+1)^{3/2} \varepsilon^{-1}.$$

Dropping the dependence on $K$ (since $K$ is fixed), this implies that $\widetilde{J}_{[\,]}(\gamma, \mathcal{P}, h) \lesssim \gamma^{1/2}$, and together with (15) we obtain $E\|\mathbb{G}_n\|_{\mathcal{M}_\gamma} \leq \sqrt{\gamma} + (\gamma\sqrt{n})^{-1}$. Since $\gamma \mapsto (\sqrt{\gamma} + (\gamma\sqrt{n})^{-1})/\gamma$ is decreasing in $\gamma$, it is a valid choice for $\phi_n(\gamma)$ in Theorem 3.4.1 of [17]. We then obtain that $r_n h(p_{\widehat{F}_n}, p_{F_0}) = O_p(1)$ provided that $h(p_{\widehat{F}_n}, p_{F_0}) \to 0$ in outer probability, and $r_n^2 \phi_n(r_n^{-1}) \le \sqrt{n}$ for all $n$. The first condition is fulfilled by the almost sure Hellinger consistency of the MLE (Theorem 3.1). The second condition holds for $r_n = cn^{1/3}$ and $c = ((\sqrt{5} - 1)/2)^{2/3}$. □

We obtain the following corollary about the $L_1(G)$ and $L_2(G)$ rates of convergence:

COROLLARY 4.2. $n^{1/3} \|\widehat{F}_n - F_0\|_{G,r} = O_p(1)$ for $r = 1, 2$.

PROOF. The result for $r = 1$ again follows from $d_{\mathrm{TV}}(p_{F_1}, p_{F_2}) \le \sqrt{2} h(p_{F_1}, p_{F_2})$. The result for $r = 2$ follows from

$$\|F - F_0\|_{G,2}^2 = \sum_{k=1}^{K+1} \int \{\sqrt{F_k} - \sqrt{F_{0k}}\}^2 \{\sqrt{F_k} + \sqrt{F_{0k}}\}^2 \, dG \le 8h^2(p_F, p_{F_0}),$$

using $\sqrt{F_k} + \sqrt{F_{0k}} \le 2$. □

4.2. *Asymptotic local minimax lower bound.* In this section we prove that $n^{1/3}$ is an asymptotic local minimax lower bound for the rate of convergence. We use the set-up of [1], Section 4.1. Let $\mathcal{P}$ be a set of probability densities on a measurable space $(\Omega, \mathcal{A})$ with respect to a $\sigma$-finite dominating measure. We estimate a parameter $\theta = Up \in \mathbb{R}$, where $U$ is a real-valued functional and $p \in \mathcal{P}$. Let $U_n$, $n \ge 1$, be a sequence of estimators based on a sample of size $n$, that is, $U_n = t_n(Z_1, \ldots, Z_n)$, where $Z_1, \ldots, Z_n$ is a sample from the density $p$, and $t_n : \Omega^n \to \mathbb{R}$ is a Borel measurable function. Let $l : [0, \infty) \to [0, \infty)$ be an increasing convex loss function with $l(0) = 0$. The risk of the estimator $U_n$ in estimating $Up$ is defined by $E_{n,p} l(|U_n - Up|)$, where $E_{n,p}$ denotes the expectation with respect to the product measure $P^{\otimes n}$ corresponding to the sample $Z_1, \ldots, Z_n$. We now recall Lemma 4.1 of [1].

LEMMA 4.3. *For any $p_1, p_2 \in \mathcal{P}$ such that the Hellinger distance $h(p_1, p_2) < 1$:*

$$\inf_{U_n} \max\{E_{n,p_1} l(|U_n - Up_1|), E_{n,p_2} l(|U_n - Up_2|)\}$$
$$\ge l(\tfrac{1}{4}|Up_1 - Up_2|(1 - h^2(p_1, p_2))^{2n}).$$

Let $k \in \{1, \ldots, K\}$ and let $U_{nk}$, $n \ge 1$, be a sequence of estimators of $F_{0k}(t_0)$. Furthermore, let $c > 0$ and let $F_n^k = (F_{n1}, \ldots, F_{nK})$ be a perturbation of $F_0$ where only the $k$th component is changed in the following way:

$$F_{nk}(x) = \begin{cases} F_{0k}(t_0 - cn^{-1/3}), & \text{if } x \in [t_0 - cn^{-1/3}, t_0), \\ F_{0k}(t_0 + cn^{-1/3}), & \text{if } x \in [t_0, t_0 + cn^{-1/3}), \\ F_{0k}(x), & \text{otherwise,} \end{cases}$$

and $F_{nj}(x) = F_{0j}(x)$ for $j \neq k$. Note that $F_n^k \in \mathcal{F}_K$ is a valid set of sub-distribution functions with overall survival function $F_{n,K+1} = 1 - F_{n+}$.

We now apply Lemma 4.3 with $l(x) = x^r$, $p_1 = p_{F_0}$ and $p_2 = p_{F_n^k}$, where $p_F$ is defined in (1). This gives a local minimax lower bound for the rate of convergence. A detailed derivation of this result is given in [10], Section 5.2.

PROPOSITION 4.4. *Fix $k \in \{1, \ldots, K\}$. Let $0 < F_{0k}(t_0) < F_{0k}(\infty)$, and let $F_{0k}$ and $G$ be continuously differentiable at $t_0$ with strictly positive derivatives $f_{0k}(t_0)$ and $g(t_0)$. Let $d = 2^{-5/3} e^{-1/3}$. Then, for $r \geq 1$,*

$$\liminf_{n \to \infty} n^{r/3} \inf_{U_n} \max\{E_{n,p_{F_0}}|U_{nk} - F_{0k}(t_0)|^r, E_{n,p_{F_n^k}}|U_{nk} - F_{nk}(t_0)|^r\}$$

$$(16) \qquad \geq d^r \left[\frac{g(t_0)}{f_{0k}(t_0)}\left\{\frac{1}{F_{0k}(t_0)} + \frac{1}{1 - F_{0+}(t_0)}\right\}\right]^{-r/3}.$$

REMARK 4.5. Note that the lower bound (16) consists of a part depending on the underlying distribution, and a universal constant $d$. It is not clear whether the constant depending on the underlying distribution is sharp, because it has not been proved that any estimator achieves this constant. However, we do know that the naive estimator $\widetilde{F}_{nk}$ does generally not achieve this constant. To see this, recall that $\widetilde{F}_{nk}$ is the MLE for the reduced data $(T_i, \Delta_k^i)$, $i = 1, \ldots, n$. Hence, its asymptotic risk is bounded below by the asymptotic local minimax lower bound for current status data:

$$d^r \left[\frac{g(t_0)}{f_{0k}(t_0)}\left\{\frac{1}{F_{0k}(t_0)} + \frac{1}{1 - F_{0k}(t_0)}\right\}\right]^{-r/3}$$

(see [1], (4.2), or take $K = 1$ in Proposition 4.4). Since $1 - F_{0k}(t_0) > 1 - F_{0+}(t_0)$ if $F_{0j}(t_0) > 0$ for some $j \in \{1, \ldots, K\}$, $j \neq k$, this bound is larger than the one given in (16).

4.3. *Local rate of convergence.* As mentioned in the introduction of this section, the $n^{1/3}$ local rate of convergence of the naive estimator and the $n^{1/3}$ local minimax lower bound for the rate of convergence suggest that the MLE converges locally at rate $n^{1/3}$. This is indeed the case, and we now give the proof of this result. However, although this result is intuitively clear, the proof is rather involved.

The two main difficulties in the proof are the lack of a closed form for the MLE and the system of sub-distribution functions. We solve the first problem by working with a characterization of the MLE in terms of necessary and sufficient conditions. This approach was also followed in [1] for case 2 interval censored data, and in [2] for convex density estimation. We handle the system of sub-distribution functions by first proving a rate result for $\widehat{F}_{n+}$ that holds uniformly on a fixed neighborhood around $t_0$, instead of on the usual shrinking neighborhood of order $n^{-1/3}$.

The outline of this section is as follows. In Section 4.3.1 we revisit the characterization of the MLE, and derive a localized version of the conditions (Proposition 4.8). In Section 4.3.2 we use this characterization to prove the rate result for $\widehat{F}_{n+}$ that is discussed above (Theorem 4.10). In Section 4.3.3 we use this result to prove the local rate of convergence for the components $\widehat{F}_{n1}, \ldots, \widehat{F}_{nK}$ (Theorem 4.17). Some technical proofs are deferred to Section 5.

Throughout, we assume that for each $k \in \{1, \ldots, K\}$, $\widehat{F}_{nk}$ is piecewise constant and right-continuous, with jumps only at points in $\mathcal{T}_k$ (see Definition 2.2). This assumption does not affect the asymptotic properties of the MLE.

4.3.1. *Revisiting the characterization.* We consider the characterization given in Corollary 2.10. Since it is difficult to work with $\widehat{F}_{nk}$ in the denominator, we start by rewriting the left-hand side of (11), using

$$\int_{[s,t)} \frac{dV_{nk}(u)}{\widehat{F}_{nk}(u)} = \int_{[s,t)} \frac{dV_{nk}(u)}{F_{0k}(u)} + \int_{[s,t)} \frac{F_{0k}(u) - \widehat{F}_{nk}(u)}{F_{0k}(u)\widehat{F}_{nk}(u)} dV_{nk}(u).$$

This leads to the following lemma:

LEMMA 4.6. *For all $k = 1, \ldots, K$ and $s, t \in \mathbb{R}$,*

$$\int_{[s,t)} \left\{ \frac{dV_{nk}(u)}{\widehat{F}_{nk}(u)} - \frac{dV_{n,K+1}(u)}{\widehat{F}_{n,K+1}(u)} \right\}$$

$$= \int_{[s,t)} \left\{ \frac{dV_{nk}(u)}{F_{0k}(u)} - \frac{dV_{n,K+1}(u)}{F_{0,K+1}(u)} \right\}$$

$$+ \int_{[s,t)} \frac{F_{0k}(u) - \widehat{F}_{nk}(u)}{F_{0k}(u)\widehat{F}_{nk}(u)} dV_{nk}(u)$$

$$- \int_{[s,t)} \frac{F_{0,K+1}(u) - \widehat{F}_{n,K+1}(u)}{F_{0,K+1}(u)\widehat{F}_{n,K+1}(u)} dV_{n,K+1}(u).$$

We now combine Corollary 2.10 and Lemma 4.6 to obtain a localized version of the characterization in Proposition 4.8. We first introduce some definitions:

DEFINITION 4.7. Let $a_k = (F_{0k}(t_0))^{-1}$ for $k = 1, \ldots, K+1$. Furthermore, for $k = 1, \ldots, K$, we define the processes $W_{nk}(\cdot)$ and $S_{nk}(\cdot)$ by

(17) $$W_{nk}(t) = \int_{u \leq t} \{\delta_k - F_{0k}(u)\} d\mathbb{P}_n(u, \delta),$$

(18) $$S_{nk}(t) = a_k W_{nk}(t) + a_{K+1} W_{n+}(t).$$

PROPOSITION 4.8. *For each $k = 1, \ldots, K$, let $0 < F_{0k}(t_0) < F_{0k}(\infty)$, and let $F_{0k}$ and $G$ be continuously differentiable at $t_0$ with strictly positive derivatives $f_{0k}(t_0)$ and $g(t_0)$. Then there is an $r > 0$ such that, for all $k = 1, \ldots, K$ and each jump point $\tau_{nk} < T_{(n)}$ of $\widehat{F}_{nk}$, we have*

$$
\begin{aligned}
(19) \quad & \int_{\tau_{nk}}^{s} \{a_k\{\widehat{F}_{nk}(u) - F_{0k}(u)\} + a_{K+1}\{\widehat{F}_{n+}(u) - F_{0+}(u)\}\} \, dG(u) \\
& \leq \int_{[\tau_{nk}, s)} dS_{nk}(u) + R_{nk}(\tau_{nk}, s) \qquad \text{for } s < T_{(n)},
\end{aligned}
$$

*where equality holds in* (19) *if $s$ is a jump point of $\widehat{F}_{nk}$, and where*

$$
(20) \quad \sup_{t_0 - 2r \leq s < t \leq t_0 + 2r} \frac{|R_{nk}(s,t)|}{n^{-2/3} \vee n^{-1/3}(t-s)^{3/2}} = O_p(1).
$$

PROOF. Let $k \in \{1, \ldots, K\}$ and let $\tau_{nk} < T_{(n)}$ be a jump point of $\widehat{F}_{nk}$. Note that Corollary 2.10 and Lemma 4.6 imply that for all $s < T_{(n)}$,

$$
\begin{aligned}
& \int_{[\tau_{nk}, s)} \frac{\widehat{F}_{nk}(u) - F_{0k}(u)}{F_{0k}(u)\widehat{F}_{nk}(u)} \, dV_{nk}(u) \\
(21) \quad & - \int_{[\tau_{nk}, s)} \frac{\widehat{F}_{n,K+1}(u) - F_{0,K+1}(u)}{F_{0,K+1}(u)\widehat{F}_{n,K+1}(u)} \, dV_{n,K+1}(u) \\
& \leq \int_{u \in [\tau_{nk}, s)} \left\{ \frac{\delta_k - F_{0k}(u)}{F_{0k}(u)} - \frac{\delta_{K+1} - F_{0,K+1}(u)}{F_{0,K+1}(u)} \right\} d\mathbb{P}_n(u, \delta),
\end{aligned}
$$

with equality if $s$ is a jump point of $\widehat{F}_{nk}$. We first consider the left-hand side of (21). For each $k \in \{1, \ldots, K+1\}$, we replace $\widehat{F}_{nk}(u)$ by $F_{0k}(u)$ in the denominator:

$$
\begin{aligned}
(22) \quad & \int_{[s,t)} \frac{\widehat{F}_{nk}(u) - F_{0k}(u)}{F_{0k}(u)\widehat{F}_{nk}(u)} \, dV_{nk}(u) \\
& = \int_{[s,t)} \frac{\widehat{F}_{nk}(u) - F_{0k}(u)}{F_{0k}(u)^2} \, dV_{nk}(u) + \rho_{nk}^{(1)}(s,t),
\end{aligned}
$$

$$
(23) \quad \text{where } \rho_{nk}^{(1)}(s,t) = -\int_{[s,t)} \frac{\{\widehat{F}_{nk}(u) - F_{0k}(u)\}^2}{F_{0k}(u)^2 \widehat{F}_{nk}(u)} \, dV_{nk}(u).
$$

Next, we replace $dV_{nk}(u)$ by $dV_k(u) = F_{0k}(u) \, dG(u)$ in the first term on the right-hand side of (22):

$$
(24) \quad \int_{[s,t)} \frac{\widehat{F}_{nk}(u) - F_{0k}(u)}{F_{0k}(u)^2} \, dV_{nk}(u) = \int_s^t \frac{\widehat{F}_{nk}(u) - F_{0k}(u)}{F_{0k}(u)} \, dG(u) + \rho_{nk}^{(2)}(s,t),
$$

$$\text{(25)} \qquad \text{where } \rho_{nk}^{(2)}(s,t) = \int_{[s,t)} \frac{\widehat{F}_{nk}(u) - F_{0k}(u)}{F_{0k}(u)^2} \, d(V_{nk} - V_k)(u).$$

Finally, we replace the denominator $F_{0k}(u)$ by $F_{0k}(t_0)$ in the first term on the right-hand side of (24):

$$\int_s^t \frac{\widehat{F}_{nk}(u) - F_{0k}(u)}{F_{0k}(u)} \, dG(u) = \int_s^t \frac{\widehat{F}_{nk}(u) - F_{0k}(u)}{F_{0k}(t_0)} \, dG(u) + \rho_{nk}^{(3)}(s,t),$$

$$\text{(26)} \qquad \text{where } \rho_{nk}^{(3)}(s,t) = \int_s^t \frac{\{\widehat{F}_{nk}(u) - F_{0k}(u)\}\{F_{0k}(t_0) - F_{0k}(u)\}}{F_{0k}(u) F_{0k}(t_0)} \, dG(u),$$

and similarly on the right-hand side of (21):

$$\int_{u \in [s,t)} \frac{\delta_k - F_{0k}(u)}{F_{0k}(u)} \, d\mathbb{P}_n(u,\delta) = \int_{u \in [s,t)} \frac{\delta_k - F_{0k}(u)}{F_{0k}(t_0)} \, d\mathbb{P}_n(u,\delta) - \rho_{nk}^{(4)}(s,t),$$

where, with $G_n$ the empirical distribution of $T_1, \ldots, T_n$ (as defined in Section 2),

$$\rho_{nk}^{(4)}(s,t) = \int_{u \in [s,t)} \frac{\{F_{0k}(u) - F_{0k}(t_0)\}\{\delta_k - F_{0k}(u)\}}{F_{0k}(u) F_{0k}(t_0)} \, d(\mathbb{P}_n - P)(u,\delta)$$

$$\text{(27)} \qquad = \int_{[s,t)} \frac{F_{0k}(u) - F_{0k}(t_0)}{F_{0k}(u) F_{0k}(t_0)} \, d(V_{nk} - V_k)(u)$$

$$\qquad + \int_{[s,t)} \frac{F_{0k}(t_0) - F_{0k}(u)}{F_{0k}(t_0)} \, d(G_n - G)(u).$$

Inequality (19) then follows from $F_{K+1} = 1 - F_+$ for $F \in \mathcal{F}_K$, and the definition

$$\text{(28)} \qquad R_{nk}(s,t) = \sum_{\ell=1}^4 \rho_{n,K+1}^{(\ell)}(s,t) - \sum_{\ell=1}^4 \rho_{nk}^{(\ell)}(s,t), \qquad k = 1, \ldots, K.$$

We now show that the remainder term $R_{nk}(s,t)$ is of the given order. Let $k \in \{1, \ldots, K+1\}$, and consider $\rho_{nk}^{(1)}$. Note that $\widehat{F}_{nk}$ and $F_{0k}$ stay away from zero with probability tending to 1 on $[t_0 - 2r, t_0 + 2r]$, by the assumption $F_{0k}(t_0) > 0$, the continuity of $F_{0k}$ at $t_0$, and the consistency of $\widehat{F}_{nk}$ (Proposition 3.3). Furthermore,

$$\int \{\widehat{F}_{nk}(u) - F_{0k}(u)\}^2 \, dV_{nk}(u)$$

$$\leq \int \{\widehat{F}_{nk}(u) - F_{0k}(u)\}^2 \, d(G_n - G)(u) + \int \{\widehat{F}_{nk}(u) - F_{0k}(u)\}^2 \, dG(u),$$

where the second term on the right-hand side is of order $O_p(n^{-2/3})$ by the $L_2(G)$ rate of convergence given in Corollary 4.2, and the first term is of order $O_p(n^{-2/3})$ by a modulus of continuity result. To see the latter, define

$$\mathcal{Q} = \{q_F(u) = \{F(u) - F_{0k}(u)\}^2 : F \in \mathcal{F}\},$$
$$\mathcal{Q}(\gamma) = \left\{q_F \in \mathcal{Q} : \int q_F(u)^2 \, dG(u) \leq \gamma^2\right\},$$

where $\mathcal{F}$ is the class of monotone functions $F : \mathbb{R} \to [0,1]$. The $L_2(G)$ rate of convergence (Corollary 4.2) implies that we can choose $C > 0$ such that $q_F \in \mathcal{Q}(Cn^{-1/3})$ with high probability. We then apply (5.42) of [16], Lemma 5.13, with $\alpha = 1$ and $\beta = 0$ to the class $\mathcal{Q}(Cn^{-1/3})$. This yields that $\rho_{nk}^{(1)}(s,t) = O_p(n^{-2/3})$ uniformly in $t_0 - 2r \leq s \leq t \leq t_0 + 2r$. Analogously, $\rho_{nk}^{(2)}(s,t) = O_p(n^{-2/3})$ uniformly in $t_0 - 2r \leq s \leq t \leq t_0 + 2r$, using the $L_2(G)$ rate of convergence and a modulus of continuity result. Next, we consider $\rho_{nk}^{(3)}(s,t)$. By the Cauchy–Schwarz inequality,

$$|\rho_{nk}^{(3)}(s,t)| \leq \left\{\int_s^t \frac{\{F_{0k}(u) - F_{0k}(t_0)\}^2}{F_{0k}(t_0)^2} \, dG(u)\right\}^{1/2}$$
$$\times \left\{\int \frac{\{\widehat{F}_{nk}(u) - F_{0k}(u)\}^2}{F_{0k}(u)^2} \, dG(u)\right\}^{1/2}.$$

The first term of the product is of order $O(t-s)^{3/2}$, uniformly in $t_0 - 2r \leq s \leq t \leq t_0 + 2r$, by the continuous differentiability of $F_{0k}$. The second term is of order $O_p(n^{-1/3})$ by the $L_2(G)$ rate of convergence. Hence, $\rho_{nk}^{(3)}(s,t) = O_p(n^{-1/3}(t-s)^{3/2})$, uniformly in $t_0 - 2r \leq s \leq t \leq t_0 + 2r$. Finally, $\rho_{nk}^{(4)}(s,t) = O_p(n^{-1/2}(t-s))$, uniformly in $t_0 - 2r \leq s \leq t \leq t_0 + 2r$, by writing $\int_{[s,t)} = \int_{[s,t_0)} - \int_{[t,t_0)}$ and using Lemma 4.9 below. Since the term $O_p(n^{-1/2}(t-s))$ is dominated by $O_p(n^{-2/3} \vee n^{-1/3}(t-s)^{3/2})$ for all $s \leq t$, it can be omitted. $\square$

LEMMA 4.9. *Let $F : \mathbb{R} \to \mathbb{R}$ be continuously differentiable at $t_0$ with derivative $f(t_0) > 0$. Then there is an $r > 0$ so that uniformly in $t_0 - 2r \leq s \leq t \leq t_0 + 2r$,*

(29) $$\int_{[s,t)} \{F(t) - F(u)\} \, d(G_n - G)(u) = O_p(n^{-1/2}(t-s)),$$

(30) $$\int_{[s,t)} \frac{F(t) - F(u)}{F(u)} \, d(V_{nk} - V_k)(u) = O_p(n^{-1/2}(t-s)),$$

$$k = 1, \ldots, K.$$

PROOF. We only prove (29), because the proof of (30) is analogous. Integration by parts yields

$$n^{1/2} \int_{[s,t)} \{F(t) - F(u)\} d(G_n - G)(u)$$
$$= -n^{1/2} \{F(t) - F(s)\} \{G_n(s) - G(s)\} + n^{1/2} \int_{[s,t)} \{G_n(u) - G(u)\} dF(u).$$

Note that $n^{1/2} \sup_{u \in \mathbb{R}} |G_n(u) - G(u)|$ is tight, since it converges in distribution to $\sup_{u \in \mathbb{R}} |B(G(u))| \leq \sup_{x \in [0,1]} |B(x)|$, where $B$ is a standard Brownian motion on $[0,1]$. Hence, both terms on the right-hand side of the display are $O_p(1)\{F(t) - F(s)\} = O_p(t-s)$, uniformly in $t_0 - 2r \leq s \leq t \leq t_0 + 2r$. □

4.3.2. *Uniform rate of convergence of $\widehat{F}_{n+}$ on a fixed neighborhood of $t_0$.* The main result of this section is a rate of convergence result for $\widehat{F}_{n+}$ which holds uniformly on a fixed neighborhood $[t_0 - r, t_0 + r]$ of $t_0$, rather than on a shrinking neighborhood of the form $[t_0 - Mn^{-1/3}, t_0 + Mn^{-1/3}]$ (Theorem 4.10). We discuss the meaning of this result in Remark 4.11, by comparing it to several existing results for current status data without competing risks. Theorem 4.10 is used in Section 4.3 to prove the local rate of convergence of the components $\widehat{F}_{n1}, \ldots, \widehat{F}_{nK}$.

THEOREM 4.10. *For all $k = 1, \ldots, K$, let $0 < F_{0k}(t_0) < F_{0k}(\infty)$, and let $F_{0k}$ and $G$ be continuously differentiable at $t_0$ with strictly positive derivatives $f_{0k}(t_0)$ and $g(t_0)$. For $\beta \in (0,1)$ we define*

$$(31) \qquad v_n(t) = \begin{cases} n^{-1/3}, & \text{if } |t| \leq n^{-1/3}, \\ n^{-(1-\beta)/3} |t|^\beta, & \text{if } |t| > n^{-1/3}. \end{cases}$$

*Then there exists a constant $r > 0$ so that*

$$(32) \qquad \sup_{t \in [t_0 - r, t_0 + r]} \frac{|\widehat{F}_{n+}(t) - F_{0+}(t)|}{v_n(t - t_0)} = O_p(1).$$

Note that the function $v_n(t) = n^{-1/3}$ for $|t| < n^{-1/3}$. Outside a $n^{-1/3}$ neighborhood we cannot expect to get a $n^{-1/3}$ rate. Therefore, for $t > n^{-1/3}$ we let the function $v_n(t)$ grow with $t$, by defining $v_n(t) = n^{-(1-\beta)/3} |t|^\beta$.

Before giving the proof of Theorem 4.10, we discuss its meaning by comparing it to several known results for current status data without competing risks.

REMARK 4.11. By taking $K = 1$ in Theorem 4.10, it follows that the theorem holds for the MLE $\widehat{F}_n$ for current status data without competing risks. Thus, to clarify the meaning of Theorem 4.10, we can compare it to known results for $\widehat{F}_n$. First, we consider the local rate of convergence given in [4], Lemma 5.4, page 95. For $M > 0$, they prove that

(33) $$\sup_{t \in [-M,M]} |\widehat{F}_n(t_0 + n^{-1/3}t) - F_0(t_0)| = O_p(n^{-1/3}).$$

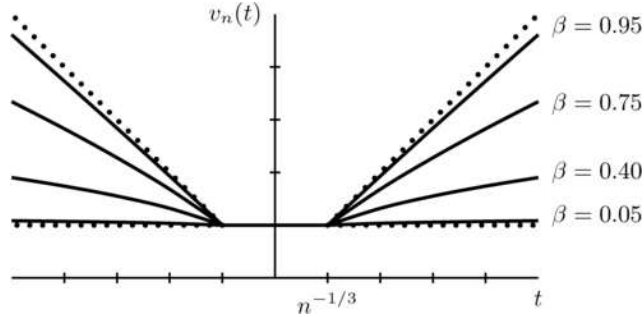

FIG. 2. *Plot of $v_n(t)$ for various values of $\beta$. The dotted lines are $y = t$ and $y = n^{-1/3}$. Note that $\beta$ close to zero gives the sharpest bound.*

We can obtain this bound by applying Theorem 4.10 to $t \in [t_0 - Mn^{-1/3}, t_0 + Mn^{-1/3}]$, and using the continuous differentiability of $F_{0k}$ at $t_0$ and the fact that

$$v_n(t - t_0) \leq v_n(Mn^{-1/3}) = M^\beta n^{-1/3}$$
$$\text{for } M \geq 1, \ t \in [t_0 - Mn^{-1/3}, t_0 + Mn^{-1/3}].$$

Hence, Theorem 4.10 implies (33) for $M \geq 1$.

Next, we consider the global bound of [4], Lemma 5.9:

(34) $$\sup_{t \in \mathbb{R}} |\widehat{F}_n(t) - F_0(t)| = O_p(n^{-1/3} \log n).$$

The result in Theorem 4.10 is fundamentally different from (34), since it is stronger than (34) for $|t - t_0| < n^{-1/3}(\log n)^{1/\beta}$, and it is weaker outside this region.

REMARK 4.12. Note that Theorem 4.10 gives a family of bounds in $\beta$. Choosing $\beta$ close to zero gives the tightest bound, as illustrated in Figure 2. For the proof of the local rate of convergence of $\widehat{F}_{n1}, \ldots, \widehat{F}_{nK}$ (Theorem 4.17), it is sufficient that Theorem 4.10 holds for one arbitrary value of $\beta \in (0,1)$. Stating the theorem for one fixed $\beta$ leads to a somewhat simpler proof. However, for completeness we present the result for all $\beta \in (0,1)$.

As an introduction to the proof of Theorem 4.10 we first note the following. Let $\varepsilon > 0$ and let $r > 0$ be small. Then the continuous differentiability of $F_{0+}$ at $t_0$ implies

$$F_{0+}(t + Mv_n(t - t_0)) \leq F_{0+}(t) + 2Mv_n(t - t_0)f_{0+}(t_0), \qquad t \in [t_0 - r, t_0 + r],$$
$$F_{0+}(t - Mv_n(t - t_0)) \geq F_{0+}(t) - 2Mv_n(t - t_0)f_{0+}(t_0), \qquad t \in [t_0 - r, t_0 + r].$$

Hence, it is sufficient to show that we can choose $n_1$ and $M$ such that for all $n > n_1$

$$P\{\exists t \in [t_0 - r, t_0 + r] : \widehat{F}_{n+}(t) \notin (F_{0+}(t - Mv_n(t - t_0)),$$
$$F_{0+}(t + Mv_n(t - t_0)))\} < \varepsilon.$$

In fact, we only prove that there exist $n_1$ and $M$ such that

(35) $\quad P\{\exists t \in [t_0, t_0 + r] : \widehat{F}_{n+}(t) \geq F_{0+}(t + Mv_n(t - t_0))\} < \varepsilon/4, \qquad n > n_1,$

since the proofs for $\widehat{F}_{n+}(t) \leq F_{0+}(t - Mv_n(t - t_0))$ and the interval $[t_0 - r, t_0]$ are analogous. In the proof of (35) we use the fact that we can choose $r$, $n_1$ and $C$ such that $P(E_{nrC}^c) < \varepsilon/8$ for all $n > n_1$, where

(36)
$$E_{nrC} = \bigcap_{k=1}^{K} \bigg\{ \widehat{F}_{nk} \text{ has a jump in } (t_0 - 2r, t_0 - r),\ T_{(n)} > t_0 + 2r,$$
$$\sup_{t_0 - 2r \leq w < t \leq t_0 + 2r} \frac{|R_{nk}(w,t)|}{n^{-2/3} \vee n^{-1/3}(t-w)^{3/2}} \leq C \bigg\},$$

and $R_{nk}(w, t)$ is defined in Proposition 4.8. For the event involving $R_{nk}$ this follows from Proposition 4.8. For the event that $\widehat{F}_{nk}$ has a jump point in $(t_0 - 2r, t_0 - r)$, this follows from consistency of $\widehat{F}_{nk}$ (Proposition 3.3) and the strict monotonicity of $F_{0k}$ in a neighborhood of $t_0$. Finally, $T_{(n)} > t_0 + 2r$ for sufficiently large $n$ follows from the positive density of $g$ in a neighborhood of $t_0$.

PROOF OF THEOREM 4.10. By the discussion above, and by writing

$$P\{\exists t \in [t_0, t_0 + r] : \widehat{F}_{n+}(t) \geq F_{0+}(t + Mv_n(t - t_0))\}$$
$$\leq P(E_{nrC}^c)$$
(37) $\qquad + P(\exists t \in [t_0, t_0 + r] : \widehat{F}_{n+}(t) \geq F_{0+}(t + Mv_n(t - t_0)), E_{nrC}),$

it is sufficient to show that we can choose $n_1$, $M$ and $C$ such that the second term of (37) is bounded by $\varepsilon/8$ for all $n > n_1$. In order to show this, we put a grid on the interval $[t_0, t_0 + r]$, analogously to [8], Lemma 4.1. The grid points $t_{nj}$ and grid cells $I_{nj}$ are denoted by

(38) $\qquad\qquad t_{nj} = t_0 + jn^{-1/3} \quad \text{and} \quad I_{nj} = [t_{nj}, t_{n,j+1})$
$$\text{for } j = 0, \ldots, J_n = \lceil rn^{1/3} \rceil.$$

This yields

$$P(\exists t \in [t_0, t_0 + r] : \widehat{F}_{n+}(t) \geq F_{0+}(t + Mv_n(t - t_0)), E_{nrC})$$
$$\leq \sum_{j=0}^{J_n} P(\exists t \in I_{nj} : \widehat{F}_{n+}(t) \geq F_{0+}(t + Mv_n(t - t_0)), E_{nrC}).$$

Hence, it is sufficient to show that we can choose $n_1$ and $m_1$ such that for all $n > n_1$, $M > m_1$ and $j = 0, \ldots, J_n$, we have

(39) $\quad P(\exists t \in I_{nj} : \widehat{F}_{n+}(t) \geq F_{0+}(t + Mv_n(t - t_0)), E_{nrC}) \leq p_{jM},$

where $p_{jM}$ satisfies $\limsup_{n \to \infty} \sum_{j=0}^{J_n} p_{jM} \to 0$ as $M \to \infty$. We prove (39) for

(40) $\quad p_{jM} = \begin{cases} d_1 \exp\{-d_2 M^3\}, & \text{if } j = 0, \\ d_1 \exp\{-d_2 (Mj^\beta)^3\}, & \text{if } j = 1, \ldots, J_n, \end{cases}$

where $d_1$ and $d_2$ are positive constants. Using the monotonicity of $\widehat{F}_{n+}$, it is sufficient to prove that for all $n > n_1$, $M > m_1$ and $j = 0, \ldots, J_n$,

(41) $\quad P\{A_{njM}, E_{nrC}\} \leq p_{jM},$

where

(42) $\quad A_{njM} = \{\widehat{F}_{n+}(t_{n,j+1}) \geq F_{0+}(s_{njM})\},$

(43) $\quad s_{njM} = t_{nj} + Mv_n(t_{nj} - t_0).$

Fix $n > 0$ and $M > 0$, and let $j \in \{0, \ldots, J_n\}$. Let $\tau_{nkj}$ be the last jump point of $\widehat{F}_{nk}$ before $t_{n,j+1}$, for $k = 1, \ldots, K$. On the event $E_{nrC}$, these jump points exist and are in $(t_0 - 2r, t_{n,j+1}]$. Without loss of generality we assume that the sub-distribution functions are labeled so that $\tau_{n1j} \leq \cdots \leq \tau_{nKj}$. On the event $A_{njM}$ there must be a $k \in \{1, \ldots, K\}$ for which $\widehat{F}_{nk}(t_{n,j+1}) \geq F_{0k}(s_{njM})$. Hence, we can define $\ell \in \{1, \ldots, K\}$ such that

(44) $\quad \widehat{F}_{nk}(t_{n,j+1}) < F_{0k}(s_{njM}), \qquad k = \ell + 1, \ldots, K,$

(45) $\quad \widehat{F}_{n\ell}(t_{n,j+1}) \geq F_{0\ell}(s_{njM}).$

Since $s_{njM} < t_0 + 2r$ for $n$ large, and $t_0 + 2r < T_{(n)}$ on the event $E_{nrC}$, we have

$$\int_{\tau_{n\ell j}}^{s_{njM}} \{a_\ell \{\widehat{F}_{n\ell}(u) - F_{0\ell}(u)\} + a_{K+1}\{\widehat{F}_{n+}(u) - F_{0+}(u)\}\} dG(u)$$
$$\leq \int_{[\tau_{n\ell j}, s_{njM})} dS_{n\ell}(u) + R_{n\ell}(\tau_{n\ell j}, s_{njM}),$$

by Proposition 4.8. Hence, $P(A_{njM}, E_{nrC})$ equals

$$P\bigg(\int_{\tau_{n\ell j}}^{s_{njM}} \{a_\ell\{\widehat{F}_{n\ell}(u) - F_{0\ell}(u)\} + a_{K+1}\{\widehat{F}_{n+}(u) - F_{0+}(u)\}\} dG(u)$$

$$\leq \int_{[\tau_{n\ell j}, s_{njM})} dS_{n\ell}(u) + R_{n\ell}(\tau_{n\ell j}, s_{njM}), A_{njM}, E_{nrC}\bigg),$$

and this is bounded above by

(46)
$$P\bigg(\int_{\tau_{n\ell j}}^{s_{njM}} a_\ell\{\widehat{F}_{n\ell}(u) - F_{0\ell}(u)\} dG(u) - \int_{[\tau_{n\ell j}, s_{njM})} dS_{n\ell}(u)$$

$$\leq R_{n\ell}(\tau_{n\ell j}, s_{njM}), A_{njM}, E_{nrC}\bigg)$$

(47)
$$+ P\bigg(\int_{\tau_{n\ell j}}^{s_{njM}} \{\widehat{F}_{n+}(u) - F_{0+}(u)\} dG(u) \leq 0, A_{njM}, E_{nrC}\bigg).$$

We now show that both terms (46) and (47) are bounded above by $p_{jM}/2$. Note that (45) implies that on the event $A_{njM}$,

$$\widehat{F}_{n\ell}(u) \geq \widehat{F}_{n\ell}(\tau_{n\ell j}) = \widehat{F}_{n\ell}(t_{n,j+1}) \geq F_{0\ell}(s_{njM}) \qquad \text{for } u \geq \tau_{n\ell j},$$

using the definition of $\tau_{n\ell j}$, and the fact that $\widehat{F}_{n\ell}$ is piecewise constant and monotone nondecreasing. Hence, on the event $A_{njM}$ we have

$$\int_{\tau_{n\ell j}}^{s_{njM}} \{\widehat{F}_{n\ell}(u) - F_{0\ell}(u)\} dG(u) \geq \int_{\tau_{n\ell j}}^{s_{njM}} \{F_{0\ell}(s_{njM}) - F_{0\ell}(u)\} dG(u)$$

$$\geq \tfrac{1}{4} g(t_0) f_{0\ell}(t_0)(s_{njM} - \tau_{n\ell j})^2,$$

for all $\tau_{n\ell j} \in [t_0 - 2r, t_{n,j+1}]$ and $r$ sufficiently small. Combining this with the definition of $E_{nrC}$ [see (36)], it follows that (46) is bounded above by

(48)
$$P\bigg(\inf_{w \in [t_0 - 2r, t_{n,j+1}]} \bigg\{\tfrac{1}{4} g(t_0) a_\ell f_{0\ell}(t_0)(s_{njM} - w)^2 - \int_{[w, s_{njM})} dS_{n\ell}(u)$$

$$- C(n^{-2/3} \vee n^{-1/3}(s_{njM} - w)^{3/2})\bigg\} \leq 0\bigg).$$

For $m_1$ and $n_1$ sufficiently large, this probability is bounded above by $p_{jM}/2$ for all $M > m_1$, $n > n_1$ and $j \in \{0, \ldots, J_n\}$, using Lemma 4.13 below. Similarly, (47) is bounded above by $p_{jM}/2$, using Lemma 4.14 below. This proves (41) and completes the proof. $\square$

Lemmas 4.13 and 4.14 play a crucial role in the proof of Theorem 4.10. The probability statement in Lemma 4.13 consists of three terms: a deterministic parabolic drift $b(s_{njM} - w)^2$, a martingale $S_{nk}$, and a remainder term $C(n^{-2/3} \vee n^{-1/6}(s_{njM} - w)^{3/2})$. The basic idea of the lemma is that the quadratic drift dominates the martingale and the remainder term. Lemma 4.14 controls the term that involves the sum of the components. In this lemma the key idea is to exploit the system of sub-distribution functions, and play out the different components against each other. The proofs of both lemmas are given in Section 5.

Finally, we note that (48) in the proof of Theorem 4.10 contains a smaller remainder term $C(n^{-2/3} \vee n^{-1/3}(s_{njM} - w)^{3/2})$ than the one in Lemma 4.13. Hence, (48) is also bounded above by $p_{jM}$. We choose to state Lemma 4.13 in terms of the larger remainder term $C(n^{-2/3} \vee n^{-1/6}(s_{njM} - w)^{3/2})$, since we need the lemma in this form for the proof of Theorem 4.17.

LEMMA 4.13. *Let $C > 0$ and $b > 0$. Then there exist $r > 0$, $n_1 > 0$ and $m_1 > 0$ such that for all $k = 1, \ldots, K$, $n > n_1$, $M > m_1$ and $j \in \{0, \ldots, J_n = \lceil rn^{1/3} \rceil\}$,*

$$P\bigg(\inf_{w \in [t_0 - 2r, t_{n,j+1}]} \bigg\{ b(s_{njM} - w)^2 - \int_{[w, s_{njM})} dS_{nk}(u) - C(n^{-2/3} \vee n^{-1/6}(s_{njM} - w)^{3/2}) \bigg\} \leq 0 \bigg) \leq p_{jM},$$

*where $s_{njM} = t_{nj} + Mv_n(t_{nj} - t_0)$, and $S_{nk}(\cdot)$, $v_n(\cdot)$ and $p_{jM}$ are defined by (18), (31) and (40), respectively.*

LEMMA 4.14. *Let the conditions of Theorem 4.10 be satisfied, and let $\ell$ be defined by (44) and (45). Then there exist $r > 0$, $n_1 > 0$ and $m_1 > 0$ such that for all $n > n_1$, $M > m_1$ and $j \in \{0, \ldots, J_n = \lceil rn^{1/3} \rceil\}$,*

$$P\bigg\{ \int_{\tau_{n\ell j}}^{s_{njM}} \{\widehat{F}_{n+}(u) - F_{0+}(u)\} dG(u) \leq 0, A_{njM}, E_{nrC} \bigg\} \leq p_{jM},$$

*where $\tau_{n\ell j}$ is the last jump point of $\widehat{F}_{n\ell}$ before $t_{n,j+1}$, $s_{njM} = t_{nj} + Mv_n(t_{nj} - t_0)$, and $E_{nrC}$, $p_{jM}$ and $A_{njM}$ are defined by (36), (40) and (42), respectively.*

REMARK 4.15. The conditions of Theorem 4.10 also hold when $t_0$ is replaced by $s$, for $s$ in a neighborhood of $t_0$. Hence, the results in this section continue to hold when $t_0$ is replaced by $s \in [t_0 - r, t_0 + r]$, for $r > 0$ sufficiently small. To be precise, there exists an $r > 0$ such that for every $\varepsilon > 0$ there exist $C > 0$ and $n_1 > 0$ such that

$$P\bigg(\sup_{t \in [t_0 - r, t_0 + r]} \frac{|\widehat{F}_{n+}(t) - F_{0+}(t)|}{v_n(t - s)} > C \bigg) < \varepsilon$$

for $s \in [t_0 - r, t_0 + r]$, $n > n_1$.

In Remark 4.12 we already mentioned that, in order to prove the local rate of convergence of the components $\widehat{F}_{n1}, \ldots, \widehat{F}_{nK}$, we only need Theorem 4.10 to hold for one value of $\beta \in (0,1)$. Therefore, we now fix $\beta = 1/2$ so that $v_n(t) = n^{-1/3} \vee n^{-1/6}\sqrt{|t|}$.

Then Remark 4.15 leads to the following corollary:

COROLLARY 4.16. *Let the conditions of Theorem 4.10 be satisfied. Then there exists an $r > 0$ such that for every $\varepsilon > 0$ there exist $C > 0$ and $n_1 > 0$ such that*

$$P\left(\sup_{t \in [t_0 - r, s]} \frac{|\int_t^s \{\widehat{F}_{n+}(u) - F_{0+}(u)\}\, dG(u)|}{n^{-2/3} \vee n^{-1/6}(s-t)^{3/2}} > C\right) < \varepsilon$$

*for $s \in [t_0 - r, t_0 + r]$, $n > n_1$.*

4.3.3. *Local rate of convergence of $\widehat{F}_{n1}, \ldots, \widehat{F}_{nK}$.* We are now ready to prove the local rate of convergence of $\widehat{F}_{n1}, \ldots, \widehat{F}_{nK}$. The proof is again based on the localized characterization given in Proposition 4.8, but we now use Corollary 4.16 to bound the term involving $\widehat{F}_{n+}$ [see (52) ahead].

THEOREM 4.17. *Let the conditions of Theorem 4.10 be satisfied. Then there exists an $r > 0$ such that for every $\varepsilon > 0$ and $M_1 > 0$ there exist $M > 0$ and $n_1 > 0$ such that*

$$P\left(\sup_{t \in [-M_1, M_1]} n^{1/3} |\widehat{F}_{nk}(s + n^{-1/3}t) - F_{0k}(s)| > M\right) < \varepsilon, \qquad k = 1, \ldots, K,$$

*for all $n > n_1$ and $s \in [t_0 - r, t_0 + r]$.*

PROOF. For the reasons discussed in Remark 4.15, it is sufficient to prove the result for $s = t_0$. Let $\varepsilon > 0$, $M_1 > 0$ and $k \in \{1, \ldots, K\}$. We want to show that there exist constants $M > M_1$ and $n_1 > 0$ such that for all $n > n_1$,

(49) $\qquad P(\widehat{F}_{nk}(t_0 + Mn^{-1/3}) \geq F_{0k}(t_0 + 2Mn^{-1/3})) < \varepsilon,$

(50) $\qquad P(\widehat{F}_{nk}(t_0 - Mn^{-1/3}) \leq F_{0k}(t_0 - 2Mn^{-1/3})) < \varepsilon.$

We only prove (49), since the proof of (50) is analogous. Define

$$B_{nkM} = \{\widehat{F}_{nk}(t_0 + Mn^{-1/3}) \geq F_{0k}(s_{nM})\} \quad \text{and} \quad s_{nM} = t_0 + 2Mn^{-1/3},$$

and let $\tau_{nk}$ be the last jump point of $\widehat{F}_{nk}$ before $t_0 + Mn^{-1/3}$. Since we may assume that $s_{nM} < t_0 + r < T_{(n)}$ for $n$ sufficiently large, Proposition 4.8 yields

$$P(B_{nkM}) = P\left(\int_{\tau_{nk}}^{s_{nM}} \{a_k \{\widehat{F}_{nk}(u) - F_{0k}(u)\}\right.$$

(51)
$$+ a_{K+1}\{\widehat{F}_{n+}(u) - F_{0+}(u)\}\} dG(u)$$
$$\leq \int_{[\tau_{nk}, s_{nM})} dS_{nk}(u) + R_{nk}(\tau_{nk}, s_{nM}), B_{nkM}\bigg).$$

By consistency of $\widehat{F}_{nk}$ (Proposition 3.3) and the strict monotonicity of $F_{0k}$ in a neighborhood of $t_0$, we may assume that $\tau_{nk} \in [t_0 - r, t_0 + Mn^{-1/3}]$. Moreover, by Proposition 4.8 and Corollary 4.16 we can choose $C > 0$ such that, with high probability,

$$|R_{nk}(\tau_{nk}, s_{nM})| \leq C(n^{-2/3} \vee n^{-1/3}(s_{nM} - \tau_{nk})^{3/2}),$$

(52) $$\int_{\tau_{nk}}^{s_{nM}} |\widehat{F}_{n+}(u) - F_{0+}(u)| dG(u) \leq C(n^{-2/3} \vee n^{-1/6}(s_{nM} - \tau_{nk})^{3/2}),$$

uniformly in $\tau_{nk} \in [t_0 - r, t_0 + Mn^{-1/3}]$. Finally, note that on the event $B_{nkM}$, we have $\int_{\tau_{nk}}^{s_{nM}} \{\widehat{F}_{nk}(u) - F_{0k}(u)\} dG(u) \geq \int_{\tau_{nk}}^{s_{nM}} \{F_{0k}(s_{nM}) - F_{0k}(u)\} dG(u)$, yielding a positive quadratic drift. The statement now follows by combining these facts with (51), and applying Lemma 4.13. □

REMARK 4.18. Note that Theorem 4.10 and Corollary 4.16 yielded the bound (52) in the proof of Theorem 4.17. Such a bound would not have been possible using rate results like (33) or (34) for $\widehat{F}_{n+}$. A bound of the form (33) cannot be used, since we cannot assume that $\tau_{nk} - s_{nM} = O_p(n^{-1/3})$. A bound of the form (34) would change the right-hand side of (52) to $Cn^{-1/3}(\tau_{nk} - s_{nM}) \log n$, and this is not dominated by the quadratic drift $(\tau_{nk} - s)^2$ for $\tau_{nk} - s > Mn^{-1/3}$. Even a stronger global bound of the form $O_p(n^{-1/3} \log \log n)$ would not suffice for this purpose. This shows that the rate result given in Theorem 4.10 was essential for the proof of Theorem 4.17.

COROLLARY 4.19. *Let the conditions of Theorem 4.10 be satisfied. For all $k = 1, \ldots, K$, let $\tau_{nk}^-(s)$ and $\tau_{nk}^+(s)$ be, respectively, the largest jump point $\leq s$ and the smallest jump point $> s$ of $\widehat{F}_{nk}$. Then there exists an $r > 0$ such that for every $\varepsilon > 0$ there exist $n_1 > 0$ and $C > 0$ such that for all $k = 1, \ldots, K$,*

$$P(\tau_{nk}^+(s) - \tau_{nk}^-(s) > Cn^{-1/3}) < \varepsilon \qquad \text{for } n > n_1, s \in [t_0 - r/2, t_0 + r/2].$$

PROOF. Let $\varepsilon > 0$ and $r > 0$. Take an arbitrary value for $M_1$ (say $M_1 = 1$), and choose $M$ and $n_1$ according to Theorem 4.17. Next, choose $C > 0$ such that

(53) $$F_{0k}(s - Cn^{-1/3}) + Mn^{-1/3} < F_{0k}(s) - Mn^{-1/3}$$
$$\text{for } s \in [t_0 - r/2, t_0 + r/2].$$

Note that $s - Cn^{-1/3} \in [t_0 - r, t_0 + r]$ for all $s \in [t_0 - r/2, t_0 + r/2]$ and $n > n_1$, for $n_1$ sufficiently large. Hence, applying Theorem 4.17 to $s$ and $s - Cn^{-1/3}$ yields

$$P(\widehat{F}_{nk}(s - Cn^{-1/3}) < F_{0k}(s - Cn^{-1/3}) + Mn^{-1/3}) > 1 - \varepsilon,$$

$$P(\widehat{F}_{nk}(s) > F_{0k}(s) - Mn^{-1/3}) > 1 - \varepsilon,$$

for $n > n_1$. Together with (53) this implies that $P(s - \tau_{nk}^-(s) > Cn^{-1/3}) < 2\varepsilon$, for $n > n_1$ and $s \in [t_0 - r/2, t_0 + r/2]$. Similar reasoning holds for $\tau_{nk}^+(s)$. □

We now obtain a bound for the remainder terms $R_{nk}(s,t)$ in Proposition 4.8, for $t_0 - mn^{-1/3} \le s \le t \le t_0 + mn^{-1/3}$ and $m > 0$. This bound is used in Proposition 3.2 of [3], which is a recentered and rescaled characterization of the MLE that is needed to prove the limiting distribution.

COROLLARY 4.20. *Let $m > 0$ and let $R_{nk}(s,t)$, $k = 1, \ldots, K$, be the remainder terms in Proposition 4.8, defined by (28). Then*

(54) $$\sup_{t_0 - mn^{-1/3} \le s \le t \le t_0 + mn^{-1/3}} |R_{nk}(s,t)| = o_p(n^{-2/3}).$$

PROOF. Since $R_{nk}(s,t) = \sum_{\ell=1}^{4} \rho_{n,K+1}^{(\ell)}(s,t) - \sum_{\ell=1}^{4} \rho_{nk}^{(\ell)}(s,t)$, it is sufficient to show that the terms $\rho_{nk}^{(\ell)}(s,t)$, $k = 1, \ldots, K+1$, $\ell = 1, \ldots, 4$, are of the right order, uniformly in $t_0 - mn^{-1/3} \le s \le t \le t_0 + mn^{-1/3}$.

Let $m > 0$ and $k \in \{1, \ldots, K+1\}$. We first consider $\rho_{nk}^{(1)}$, defined by (23). By the local rate of convergence (Theorem 4.17) and the continuous differentiability of $F_{0k}$ at $t_0$, we have $\widehat{F}_{nk}(u) - F_{0k}(u) = O_p(n^{-1/3})$, uniformly in $u \in [t_0 - mn^{-1/3}, t_0 + mn^{-1/3}]$. Moreover, the assumption $F_{0k}(t_0) > 0$, the consistency of $\widehat{F}_{nk}$ (Proposition 3.3), and the continuity of $F_{0k}$ at $t_0$, imply that $\{F_{0k}(u)\widehat{F}_{nk}(u)\}^{-1} = O_p(1)$, uniformly in $u \in [t_0 - mn^{-1/3}, t_0 + mn^{-1/3}]$. Hence,

$$|\rho_{nk}^{(1)}(s,t)| \le O_p(n^{-2/3}) \int_{[t_0 - mn^{-1/3}, t_0 + mn^{-1/3})} dV_{nk}(u) = O_p(n^{-1}),$$

uniformly in $t_0 - mn^{-1/3} \le s \le t \le t_0 + mn^{-1/3}$.

Next, we consider $\rho_{nk}^{(2)}$, defined by (25). We apply Theorem 2.11.22 of [17] to the class $\mathcal{Q}_n$, where

$$\mathcal{Q}_n = \left\{ q_{n, F_n, t}(u) = \sqrt{n} \frac{F_n(u) - F_{0k}(u)}{F_{0k}(u)^2} 1_{[t_0, t_0 + n^{-1/3}t)}(u) : t \in [-m, m], F_n \in \mathcal{F}_n \right\},$$

$$\mathcal{F}_n = \Big\{ F_n : \mathbb{R} \mapsto [0,1],$$

$$F_n \text{ monotone}, \sup_{u \in [-m, m]} |(F_n - F_{0k})(t_0 + n^{-1/3}u)| \le Cn^{-1/3} \Big\}.$$

This yields that the sequence $\{\sqrt{n}(V_{nk} - V_k)q_{n,F_n,t} : t \in [-m,m], F_n \in \mathcal{F}_n\}$ is tight. Moreover, for every $\varepsilon > 0$ we can choose $C > 0$ and $n_1 > 0$ such that $P(\widehat{F}_{nk} \in \mathcal{F}_n) > 1 - \varepsilon$ for all $n > n_1$, by the local rate of convergence of $\widehat{F}_{nk}$ (Theorem 4.17) and the continuous differentiability of $F_{0k}$ at $t_0$. This implies that $\rho_{nk}^{(2)}(s,t) = O_p(n^{-1})$, uniformly in $t_0 - mn^{-1/3} \leq s \leq t \leq t_0 + mn^{-1/3}$, since

$$\sqrt{n}(V_{nk} - V_k)q_{n,\widehat{F}_{nk},t} = n \int_{[t_0, t_0 + n^{-1/3}t)} \frac{\widehat{F}_{nk}(u) - F_{0k}(u)}{F_{0k}(u)^2} \, d(V_{nk} - V_k)(u).$$

Finally, we consider the terms $\rho_{nk}^{(3)}$ and $\rho_{nk}^{(4)}$, defined by (26) and (27). We showed in the proof of Proposition 4.8 that $\rho_{nk}^{(3)}(s,t) = O_p(n^{-1/3}(t-s)^{3/2})$ and $\rho_{nk}^{(4)}(s,t) = O_p(n^{-1/2}(t-s))$, uniformly in $t_0 - r \leq s \leq t \leq t_0 + r$. Plugging in $t - s < 2mn^{-1/3}$ completes the proof. □

## 5. Technical proofs.

PROOF OF LEMMA 4.13. Let $k \in \{1, \ldots, K\}$, $n > 0$ and $j \in \{0, \ldots, J_n\}$. Note that for $M$ large, we have for all $w \leq t_{n,j+1}$:

$$C(n^{-2/3} \vee n^{-1/6}(s_{njM} - w)^{3/2}) \leq \tfrac{1}{2}b(s_{njM} - w)^2,$$

since $s_{njM} - w \geq (M-1)n^{-1/3}$. Hence, the probability in the statement of Lemma 4.13 is bounded above by

$$(55) \quad P\left\{\sup_{w \in [t_0 - 2r, t_{n,j+1}]} \left\{\int_{[w, s_{njM})} dS_{nk}(u) - \tfrac{1}{2}b(s_{njM} - w)^2\right\} \geq 0\right\}.$$

In order to bound this probability, we put a grid on the interval $[t_0 - 2r, t_{n,j+1})$, with grid points $t_{n,j-q}$ and grid cells $I_{n,j-q}$ given by

$$(56) \quad \begin{aligned} I_{n,j-q} &= [t_{n,j-q}, t_{n,j-q+1}) \\ &= [t_0 + (j-q)n^{-1/3}, t_0 + (j-q+1)n^{-1/3}), \end{aligned}$$

for $q = 0, \ldots, Q_{nj} = \lceil 2rn^{1/3} + j \rceil$. Then (55) is bounded above by

$$(57) \quad \sum_{q=0}^{Q_{nj}} P\left\{\sup_{w \in I_{n,j-q}} \int_{[w, s_{njM})} dS_{nk}(u) \geq \tfrac{1}{2}b(s_{njM} - t_{n,j-q+1})^2\right\}.$$

If we bound the $q$th term in (57) by

$$(58) \quad p_{jqM} = \begin{cases} \exp\{-d_2(q+M)^3\}, & \text{if } j = 0,\ q = 0, \ldots, Q_{n0}, \\ \exp\{-d_2(q+Mj^\beta)^3\}, & \text{if } j = 1, \ldots, J_n,\ q = 0, \ldots, Q_{nj}, \end{cases}$$

for some $d_2 > 0$, then we are done, since summing over $q$ and using $(a+b)^3 \geq a^3 + b^3$ for $a, b > 0$, and defining $d_1 = \sum_{q=0}^{\infty} \exp(-d_2 q^3) < \infty$, yields

$$p_{jM} \leq \begin{cases} d_1 \exp\{-d_2 M^3\}, & \text{if } j = 0, \\ d_1 \exp\{-d_2 (Mj^\beta)^3\}, & \text{if } j = 1, \ldots, J_n. \end{cases}$$

In order to prove that such a bound holds, we introduce, for each $\theta > 0$, the time-reversed submartingale $\exp\{n\theta \int_{[w, s_{njM})} dS_{nk}(u)\}$, for $w \leq s_{njM}$, with respect to the filtration $\{\mathcal{F}_w : w \leq t_0 + r\}$, where $\mathcal{F}_w = \sigma\{(T_i, \Delta^i), i = 1, \ldots, n : T_i \geq w\}$. Then, by Doob's submartingale inequality (see, e.g., [12], Theorem 70.1, page 177), the $q$th term in (57) is, for each $\theta > 0$, bounded above by

$$P\left\{\sup_{w \in I_{n,j-q}} \exp\left\{n\theta \int_{[w, s_{njM})} dS_{nk}(u)\right\} \geq \exp\{\tfrac{1}{2} n\theta b(s_{njM} - t_{n,j-q+1})^2\}\right\}$$

$$(59) \quad \leq \exp\{-\tfrac{1}{2} n\theta b(s_{njM} - t_{n,j-q+1})^2\} E \exp\left\{n\theta \int_{[t_{n,j-q}, s_{njM})} dS_{nk}(u)\right\}.$$

We are now left with computing an upper bound for $E \exp\{n\theta \int_{[t_{n,j-q}, s_{njM})} dS_{nk}(u)\}$. Since we have i.i.d. observations, this expectation can be written as

$$(60) \quad (E \exp\{\theta 1_{[t_{n,j-q}, s_{njM})}(T) \zeta_{nk}(T, \Delta)\})^n$$

$$\text{where } \zeta_{nk}(T, \Delta) = \frac{\Delta_k}{F_{0k}(T)} - \frac{\Delta_{K+1}}{F_{0,K+1}(T)}.$$

Using the exponential series and $E(\zeta_{nk}(T, \Delta)|T) = 0$, (60) equals

$$\exp\left\{n \log\left(1 + E 1_{[t_{n,j-q}, s_{njM})}(T) \sum_{\ell=2}^{\infty} \frac{\theta^\ell \zeta_{nk}(T, \Delta)^\ell}{\ell!}\right)\right\},$$

and since $\log(1 + x) \leq x$ for all $x > -1$, this is bounded above by

$$(61) \quad \exp\{\tfrac{1}{2} n f_n(\theta, t_{n,j-q}, s_{njM}) \theta^2 (s_{njM} - t_{n,j-q})\}$$

$$\text{where } f_n(\theta, c_1, c_2) \equiv \frac{2}{c_2 - c_1} \sum_{\ell=2}^{\infty} \frac{\theta^{\ell-2}}{\ell!} \int_{c_1}^{c_2} |E\{\zeta_{nk}(T, \Delta)^\ell | T = t\}| \, dG(t).$$

Next, for each pair $c_1 < c_2$, we let $\theta_{c_1, c_2}$ be the solution of the equation $\theta f_n(\theta, c_1, c_2) = \tfrac{1}{4} b(c_2 - c_1)$. This solution exists and is unique for all $c_1 < c_2$, since $\theta \mapsto \theta f_n(\theta, c_1, c_2)$ is a continuous increasing map from $\mathbb{R}_+$ onto $\mathbb{R}_+$. Choosing $\theta = \theta_{t_{n,j-q}, s_{njM}}$ in (61), and using that $(s_{njM} - t_{n,j-q})^2 \leq 2(s_{njM} - t_{n,j-q+1})^2$ for all $j$ and $q$ and $M > 4$, and that $s_{njM} - t_{n,j-q} \geq s_{njM} - t_{n,j-q+1}$, yields that (59) is bounded above by

$$(62) \quad \exp\left\{-\frac{nb^2(s_{njM} - t_{n,j-q+1})^3}{16 f_n(\theta_{t_{n,j-q}, s_{njM}}, t_{n,j-q}, s_{njM})}\right\} \leq \exp\left\{-\frac{nb^2(s_{njM} - t_{n,j-q+1})^3}{16d}\right\},$$

where $d \equiv \sup_{t_0-2r \leq c_1 < c_2 \leq t_0+2r} f_n(\theta_{c_1,c_2}, c_1, c_2)$. Here we use that, for $n$ sufficiently large, all intervals $[t_{n,j-q}, s_{njM})$ are contained in the interval $[t_0 - 2r, t_0 + 2r]$. Note that $d < \infty$ since

$$\theta_{c_1,c_2} \leq \tfrac{1}{4} b(c_2-c_1)^2 \Big/ \int_{c_1}^{c_2} |E\{\zeta_{nk}(T,\Delta)^2 | T=t\}| \, dG(t).$$

Hence, there is a constant $d_2 > 0$ such that for all $q = 0, \ldots, Q_{nj}$, the right-hand side of (62) is bounded above by $\exp\{-d_2(q+M)^3\}$ for $j = 0$, and by $\exp\{-d_2(q+Mj^\beta)^3\}$ for $j = 1, \ldots, J_n$. $\square$

PROOF OF LEMMA 4.14. We first note that $\ell$ is only defined on the event $A_{njM} = \{\widehat{F}_{n+}(t_{n,j+1}) \geq F_{0+}(s_{njM})\}$. Hence, this entire proof should be read on the event $A_{njM}$. Furthermore, note that the lemma is trivial if $\ell = K$, because in that case $\widehat{F}_{n+}(u) \geq F_{0+}(s_{njM})$ for all $u \geq \tau_{n\ell j}$. Therefore, suppose $\ell < K$. Then we typically do not have that $\widehat{F}_{n+}(u) \geq F_{0+}(s_{njM})$ for all $u \geq \tau_{n\ell j}$, since $\widehat{F}_{n+}(u)$ may have jumps on $(\tau_{n\ell j}, t_{n,j+1})$. We now exploit the $K$-dimensional system of sub-distribution functions by breaking $\int_{\tau_{n\ell j}}^{s_{njM}} \{\widehat{F}_{n+}(u) - F_{0+}(u)\} \, dG(u)$ into pieces that we analyze separately. First, we define $\ell^* \in \{\ell, \ldots, K\}$ as follows. If

(63) $$\int_{\tau_{n\ell j}}^{\tau_{nkj}} \{\widehat{F}_{n+}(u) - F_{0+}(u)\} \, dG(u) \leq 0 \quad \text{for all } k = \ell+1, \ldots, K,$$

we let $\ell^* = \ell$. Otherwise we define $\ell^*$ such that

(64) $$\int_{\tau_{n\ell j}}^{\tau_{nkj}} \{\widehat{F}_{n+}(u) - F_{0+}(u)\} \, dG(u) \leq 0, \quad k = \ell^*+1, \ldots, K,$$

(65) $$\int_{\tau_{n\ell j}}^{\tau_{n\ell^* j}} \{\widehat{F}_{n+}(u) - F_{0+}(u)\} \, dG(u) > 0.$$

Then, by (65) and the decomposition $\int_{\tau_{n\ell j}}^{s_{njM}} = \int_{\tau_{n\ell j}}^{\tau_{n\ell^* j}} + \int_{\tau_{n\ell^* j}}^{s_{njM}}$, we get

(66) $$\int_{\tau_{n\ell j}}^{s_{njM}} \{\widehat{F}_{n+}(u) - F_{0+}(u)\} \, dG(u)$$
$$\geq \int_{\tau_{n\ell^* j}}^{s_{njM}} \{\widehat{F}_{n+}(u) - F_{0+}(u)\} \, dG(u),$$

where strict inequality holds if $\ell \neq \ell^*$. By rearranging the sum and using the notation $\tau_{n,K+1,j} = s_{njM}$, we can write the right-hand side of (66) as

(67) $$\sum_{k=\ell^*+1}^{K} \int_{\tau_{n\ell^* j}}^{\tau_{nkj}} \{\widehat{F}_{nk}(u) - F_{0k}(u)\} \, dG(u)$$
$$+ \sum_{k=\ell^*}^{K} \sum_{p=1}^{k} \int_{\tau_{nkj}}^{\tau_{n,k+1,j}} \{\widehat{F}_{np}(u) - F_{0p}(u)\} \, dG(u).$$

We now derive lower bounds for both terms in (67), on the event $A_{njM} \cap E_{nrC}$. Starting with the first term, note that

$$\text{(68)} \quad \int_{\tau_{n\ell^* j}}^{\tau_{nkj}} \{\widehat{F}_{n+}(u) - F_{0+}(u)\} \, dG(u) \leq 0, \qquad k = \ell^* + 1, \ldots, K.$$

Namely, if $\ell = \ell^*$, then (68) is the same as (63). On the other hand, if $\ell < \ell^*$, then (68) follows (with strict inequality) from (64), (65) and the decomposition $\int_{\tau_{n\ell j}}^{\tau_{nkj}} = \int_{\tau_{n\ell j}}^{\tau_{n\ell^* j}} + \int_{\tau_{n\ell^* j}}^{\tau_{nkj}}$. Furthermore, Proposition 4.8 implies that on the event $E_{nrC}$,

$$\text{(69)} \quad \int_t^{\tau_{nkj}} \{a_k \{\widehat{F}_{nk}(u) - F_{0k}(u)\} + a_{K+1} \{\widehat{F}_{n+}(u) - F_{0+}(u)\}\} \, dG(u)$$
$$\geq \int_{[t,\tau_{nkj})} dS_{nk}(u) - C(n^{-2/3} \vee n^{-1/3}(\tau_{nkj} - t)^{3/2}),$$

for $k = 1, \ldots, K$ and $t < \tau_{nkj}$, where $S_{nk}$ is defined in (18). Using this inequality with $t = \tau_{n\ell^* j}$ together with (68) yields that on the event $E_{nrC}$,

$$\int_{\tau_{n\ell^* j}}^{\tau_{nkj}} a_k \{\widehat{F}_{nk}(u) - F_{0k}(u)\} \, dG(u)$$
$$\geq \int_{[\tau_{n\ell^* j}, \tau_{nkj})} dS_{nk}(u) - C(n^{-2/3} \vee n^{-1/3}(\tau_{nkj} - \tau_{n\ell^* j})^{3/2}),$$

for $k = \ell^* + 1, \ldots, K$, so that the first term of (67) is bounded below by

$$\sum_{k=\ell^*+1}^{K} a_k^{-1} \left\{ \int_{[\tau_{n\ell^* j}, \tau_{nkj})} dS_{nk}(u) - C(n^{-2/3} \vee n^{-1/3}(\tau_{nkj} - \tau_{n\ell^* j})^{3/2}) \right\}.$$

We now derive a lower bound for the second term of (67). Note that the inequalities (44) in the definition of $\ell$ imply that on the event $A_{njM}$

$$\sum_{p=k+1}^{K} \widehat{F}_{np}(t_{n,j+1}) < \sum_{p=k+1}^{K} F_{0p}(s_{njM}), \qquad k = \ell, \ldots, K.$$

Together with the definition of $\tau_{n1j}, \ldots, \tau_{nKj}$, this yields that on the event $A_{njM} = \{\widehat{F}_{n+}(t_{n,j+1}) \geq F_{0+}(s_{njM})\}$, we have

$$\sum_{p=1}^{k} \widehat{F}_{np}(\tau_{npj}) = \sum_{p=1}^{k} \widehat{F}_{np}(t_{n,j+1})$$
$$> \sum_{p=1}^{k} F_{0p}(s_{njM}), \qquad k = \ell, \ldots, K.$$

Furthermore, $\widehat{F}_{np}(\tau_{npj}) \leq \widehat{F}_{np}(\tau_{nkj})$ for $p \leq k$ by the monotonicity of $\widehat{F}_{np}$ and the ordering $\tau_{n1j} \leq \ldots \leq \tau_{nKj}$. Hence, we get for $k = \ell, \cdots, K$ and $u \geq \tau_{nkj}$:

$$\sum_{p=1}^{k} \widehat{F}_{np}(u) \geq \sum_{p=1}^{k} \widehat{F}_{np}(\tau_{nkj})$$

$$\geq \sum_{p=1}^{k} \widehat{F}_{n+}(\tau_{np}) > \sum_{p=1}^{k} F_{0p}(s_{njM}).$$

This implies that the second term of (67) is bounded below by

$$\sum_{k=\ell^*}^{K} \sum_{p=1}^{k} \int_{\tau_{nkj}}^{\tau_{n,k+1,j}} \{F_{0p}(s_{njM}) - F_{0p}(u)\} dG(u)$$

$$= \sum_{k=1}^{K} \int_{\tau_{nkj} \vee \tau_{n\ell^*j}}^{s_{njM}} \{F_{0k}(s_{njM}) - F_{0k}(u)\} dG(u).$$

Hence,

$$P\left\{\int_{\tau_{n\ell j}}^{s_{njM}} \{\widehat{F}_{n+}(u) - F_{0+}(u)\} dG(u) \leq 0,\, A_{njM},\, E_{nrC}\right\}$$

$$\leq P\Bigg\{\sum_{k=\ell^*+1}^{K} a_k^{-1} \bigg\{\int_{[\tau_{n\ell^*j}, \tau_{nkj})} dS_{nk}(u)$$

$$- C(n^{-2/3} \vee n^{-1/3}(\tau_{nkj} - \tau_{n\ell^*j})^{3/2})\bigg\}$$

$$+ \sum_{k=1}^{K} \int_{\tau_{nkj} \vee \tau_{n\ell^*j}}^{s_{njM}} \{F_{0k}(s_{njM}) - F_{0k}(u)\} dG(u) \leq 0,\, E_{nrC}\Bigg\}.$$

The statement now follows by writing

$$\int_{[\tau_{n\ell^*j}, \tau_{nkj})} dS_{nk}(u) = \int_{[\tau_{n\ell^*j}, s_{njM})} dS_{nk}(u) - \int_{[\tau_{nkj}, s_{njM})} dS_{nk}(u)$$

and several applications of Lemma 4.13. $\square$

## REFERENCES

[1] GROENEBOOM, P. (1996). Lectures on inverse problems. *Lectures on Probability Theory and Statistics. Ecole d'Eté de Probabilités de Saint Flour XXIV—1994. Lecture Notes in Math.* **1648** 67–164. Springer, Berlin. MR1600884
[2] GROENEBOOM, P., JONGBLOED, G. and WELLNER, J. A. (2001). Estimation of a convex function: Characterizations and asymptotic theory. *Ann. Statist.* **29** 1653–1698. MR1891742

[3] GROENEBOOM, P., MAATHUIS, M. H. and WELLNER, J. A. (2008). Current status data with competing risks: Limiting distribution of the MLE. *Ann. Statist.* **36** 1064–1089.

[4] GROENEBOOM, P. and WELLNER, J. A. (1992). *Information Bounds and Nonparametric Maximum Likelihood Estimation*. Birkhäuser, Basel. MR1180321

[5] HUDGENS, M. G., SATTEN, G. A. and LONGINI, JR., I. M. (2001). Nonparametric maximum likelihood estimation for competing risks survival data subject to interval censoring and truncation. *Biometrics* **57** 74–80. MR1821337

[6] JEWELL, N. P. and VAN DER LAAN, M. J. (2004). Current status data: Review, recent developments and open problems. In *Advances in Survival Analysis. Handbook of Statist.* **23** 625–642. North-Holland, Amsterdam. MR2065792

[7] JEWELL, N. P., VAN DER LAAN, M. J. and HENNEMAN, T. (2003). Nonparametric estimation from current status data with competing risks. *Biometrika* **90** 183–197. MR1966559

[8] KIM, J. and POLLARD, D. (1990). Cube root asymptotics. *Ann. Statist.* **18** 191–219. MR1041391

[9] MAATHUIS, M. H. (2005). Reduction algorithm for the MLE for the distribution function of bivariate interval censored data. *J. Comput. Graph. Statist.* **14** 352–362. MR2160818

[10] MAATHUIS, M. H. (2006). Nonparametric estimation for current status data with competing risks. Ph.D. thesis, Univ. Washington. Available at http://stat.ethz.ch/~maathuis/papers/.

[11] PFANZAGL, J. (1988). Consistency of maximum likelihood estimators for certain nonparametric families, in particular: Mixtures. *J. Statist. Plann. Inference* **19** 137–158. MR0944202

[12] ROGERS, L. C. G. and WILLIAMS, D. (1994). *Diffusions, Markov Processes, and Martingales.* **1**, 2nd ed. Wiley, Chichester. MR1331599

[13] SCHICK, A. and YU, Q. (2000). Consistency of the GMLE with mixed case interval-censored data. *Scand. J. Statist.* **27** 45–55. MR1774042

[14] VAN DE GEER, S. A. (1993). Hellinger-consistency of certain nonparametric maximum likelihood estimators. *Ann. Statist.* **21** 14–44. MR1212164

[15] VAN DE GEER, S. A. (1996). Rates of convergence of the maximum likelihood estimator in mixture models. *J. Nonparametr. Statist.* **6** 293–310. MR1386341

[16] VAN DE GEER, S. A. (2000). *Applications of Empirical Process Theory.* Cambridge Univ. Press. MR1739079

[17] VAN DER VAART, A. W. and WELLNER, J. A. (1996). *Weak Convergence and Empirical Processes: With Applications to Statistics.* Springer, New York. MR1385671

[18] VAN DER VAART, A. W. and WELLNER, J. A. (2000). Preservation theorems for Glivenko–Cantelli and uniform Glivenko–Cantelli classes. In *High Dimensional Probability II* 115–133. Birkhäuser, Boston. MR1857319

[19] ZEIDLER, E. (1985). *Nonlinear Functional Analysis and its Applications III: Variational Methods and Optimization.* Springer, New York. MR0768749

P. GROENEBOOM
DEPARTMENT OF MATHEMATICS
DELFT UNIVERSITY OF TECHNOLOGY
MEKELWEG 4
2628 CD DELFT
THE NETHERLANDS
E-MAIL: p.groeneboom@ewi.tudelft.nl

M. H. MAATHUIS
SEMINAR FÜR STATISTIK
ETH ZÜRICH
CH-8092 ZÜRICH
SWITZERLAND
USA
E-MAIL: maathuis@stat.math.ethz.ch

J. A. Wellner
Department of Statistics
University of Washington
Seattle, Washington 98195
USA
E-mail: jaw@stat.washington.edu